\documentclass[onefignum,onetabnum]{siamart251216}

\usepackage{amssymb,mathtools,bm,bbm,mathrsfs}

\usepackage{tabularx,array,multirow,booktabs,diagbox}
\usepackage{float}
\usepackage[section]{placeins}
\usepackage{multicol}
\usepackage{silence}
\usepackage{caption}
\usepackage{subcaption}

\usepackage{graphicx}
\graphicspath{{fig/}}
\usepackage[export]{adjustbox}
\usepackage{tikz}
\usetikzlibrary{plotmarks}

\usepackage{xcolor}
\usepackage[normalem]{ulem}
\usepackage{comment}
\usepackage{todonotes}
\usepackage{enumerate}

\usepackage{pifont}

\usepackage{siunitx}
\newsiamthm{problem}{Problem}
\newsiamremark{remark}{Remark}
\newsiamremark{example}{Example}
\usepackage{bm}
 
\newcommand{\Om}{\Omega}
\newcommand{\Omt}{\Omega_2}
\newcommand{\Omo}{\Omega_1}
\newcommand{\Hoo}{H^1_0(\Om)}
\newcommand{\Hot}{H^1(\Omt)}
\newcommand{\Lam}{\Lambda}
\newcommand{\R}{R}                 %
\newcommand{\E}{E}                 %
\newcommand{\W}{W}
\newcommand{\Q}{\Lam}
\newcommand{\inner}[2]{\left(#1,#2\right)}
\newcommand{\dual}[2]{\left\langle #1,#2\right\rangle}
\newcommand{\norm}[1]{\left\|#1\right\|}
\newcommand{\seminorm}[1]{\left|#1\right|}
\newcommand{\betaone}{\beta_1}
\newcommand{\rtwo}{\beta_2}
\newcommand{\gap}{\alpha}
\newcommand{\Bcal}{\mathcal{B}}
\newcommand{\Acal}{\mathcal{A}}
\newcommand{\bbA}{\mathbb{A}}
\newcommand{\bbB}{\mathbb{B}}
\newcommand{\bbBt}{\widetilde{\mathbb{B}}}
\newcommand{\mat}[1]{\mathsf{#1}}   %
\newcommand{\equeq}{\simeq}
\DeclareMathOperator{\diag}{diag}
\newcommand{\ind}{\mathbf{1}}
\DeclareRobustCommand{\methodkey}[2][]{%
    \tikz[baseline=-0.5ex, x=1em, y=1em]{%
        \draw[black, semithick, #1] (0,0) -- (1.8,0);
        \draw[black, semithick] plot[mark=#2, mark size=0.25em] coordinates {(0.9,0)};
    }%
}

\headers{Operator preconditioning for FD-DLM elliptic interface problems}{M.~Feder}

\title{Operator preconditioning for fictitious domain discretizations of elliptic interface problems}

\author{Marco Feder\thanks{Department of Mathematics, University of Pisa, Largo B. Pontecorvo 5, 56126 Pisa, Italy (\email{marco.feder@dm.unipi.it}). Corresponding author: M.~Feder.}\funding{MF acknowledges support
of the European Research Council (ERC) under the European Union's Horizon 2020 research and innovation programme (call HORIZON-EUROHPC-JU-2023-COE-03, grant agreement No. 101172493 "dealii-X").
The author thanks the Barcelona Supercomputing Center for the computational resources made available through the EuroHPC projects EHPC-DEV-2026D02-275 and EHPC-DEV-2026D07-202. The author is a member
of Gruppo Nazionale per il Calcolo Scientifico (GNCS), Istituto Nazionale di Alta Matematica (INdAM).
}}

\begin{document}

\maketitle

\begin{abstract}
    We develop a parameter-robust preconditioner for saddle point systems arising from fictitious domain discretizations with distributed Lagrange multipliers of elliptic interface
    problems. Using the operator preconditioning framework, we equip the function spaces with weighted norms and establish
    continuous and discrete stability with constants independent of the material parameters and their contrast.
    The resulting preconditioner is block-diagonal and consists only of standard elliptic operators, allowing
    efficient inversion through multilevel methods. Numerical experiments on non-matching meshes confirm robustness
    with respect to mesh refinement and coefficient contrast. In addition, we present a parallel implementation combining
    matrix-free operator evaluations with geometric multigrid preconditioning, demonstrating good strong scalability.
\end{abstract}

\begin{keywords}
    operator preconditioning, iterative solvers, mixed finite elements, fictitious domain methods.
\end{keywords}

\begin{MSCcodes}
    65F08, 65F10, 65N30, 65N85
\end{MSCcodes}

\section{Introduction}\label{sec:intro}
Fictitious domain methods solve partial differential equations on complex or moving geometries by
embedding the physical domain into a simpler background domain and enforcing coupling conditions
through suitable constraints~\cite{Glowinski}.
In the fictitious domain method with distributed Lagrange multipliers (FD-DLM), introduced for
elliptic interface problems in~\cite{AURICCHIO201536,BOFFI20142151}, the constraint is imposed by a
Lagrange multiplier distributed over the immersed domain, yielding a $3\times3$ block system. Inspired by the immersed boundary method~\cite{Peskin_2002}, this
formulation extends naturally to fluid--structure interaction systems~\cite{BoffiCavalliniGastaldi2015,BG_numermath}, and its
stability, approximation properties, and alternative formulations have been studied
extensively in recent years~\cite{BoffiGastaldi2024,BoffiConditioning,REGAZZONI2024117327}.

The resulting multiple saddle point systems become increasingly ill-conditioned under mesh refinement
and coefficient contrast, motivating the development of specialized solvers~\cite{BOFFI2024406,alshehri2025multigridpreconditioningfddlmmethod}. In recent joint works, augmented Lagrangian (AL)
preconditioners for fictitious domain discretizations in codimension-one~\cite{BENZI2026118522} and FD-DLM discretizations of interface problems~\cite{BENZI2026522} have been
developed, and then extended to mixed-dimensional beam--solid coupling in~\cite{FIRMBACH2027119425}. Overall, they yield low iteration counts with flexible GMRES (FGMRES), albeit with a certain
sensitivity to the choice of the augmentation parameters. This robustness, however, comes at a
price. The augmented diagonal blocks might involve triple matrix--matrix products that could modify the sparsity
pattern of the original blocks, and their (approximate) inversion requires specialized multigrid solvers~\cite{ALprec,Farrell2019}, based on
the ideas introduced in~\cite{schoeberl_1999}. In some cases,
the explicit construction of triple matrix--matrix products can be avoided entirely~\cite{benzi2026unifiedvariationalframeworkaugmented}. Nevertheless, optimal multilevel schemes for
the augmented blocks arising in fictitious domain discretizations have not yet been developed, and their design
is a nontrivial task in its own right due to the presence of overlapping meshes.
In addition, since the preconditioned systems are nonsymmetric, eigenvalue information alone does not determine the convergence of
GMRES methods~\cite{GMRES_eigs}, and a formal convergence proof would require a
dedicated field-of-values analysis.%

To overcome these challenges, we instead apply the operator preconditioning framework from~\cite{MW}. Its main ingredient is the use of
mapping properties of the underlying continuous operators to derive block-diagonal preconditioners for the systems derived from the discretization procedure. Related approaches have been recently developed for
coupled systems~\cite{hu2026parameter,cai2026parameter,sande2025robust} and for
codimension-one and codimension-two interface couplings~\cite{prec1d2d,haznics,dimola2024robust}. In the latter
settings, however, the multiplier block generally requires the realization of fractional operators. In this work, we construct a block-diagonal,
symmetric positive definite preconditioner compatible with MinRes~\cite{PaigeSaunders1975}. We introduce coefficient-dependent norms that yield
continuous and, under a mesh ratio condition, discrete stability constants independent of the material coefficients and their contrast. The diagonal blocks consist only of
standard second-order elliptic operators, avoiding fractional operators, augmentation parameters, and triple matrix--matrix products,
and are amenable to order-optimal multilevel inversion through highly optimized matrix-free techniques.
Compared with the block-triangular augmented Lagrangian preconditioner developed in~\cite{BENZI2026522}, the proposed method
requires more outer iterations but has a substantially lower cost per iteration, leading in practice to considerably shorter solution times.

To demonstrate the practical potential of the proposed approach, we perform large-scale parallel
experiments with the \textsc{deal.II} library~\cite{dealIIdesign}, exploiting state-of-the-art
matrix-free operator evaluations from~\cite{KRONBICHLER2012135} within geometric multigrid
preconditioners for the inversion of the diagonal blocks, solving three-dimensional problems with more than $400$ million degrees of freedom. To the best of the author's knowledge, this is the first time that the FD-DLM scheme has been realized in
a matrix-free framework, as well as the first time it has been tested on problems of this size.

The outline of the paper is as follows: Section~\ref{sec:model_problem} introduces the
model problem, Section~\ref{sec:operator_prec} presents the application of the operator preconditioning
framework. Section~\ref{sec:discretization} discusses the finite element discretization of the problem, the resulting linear system, the associated
preconditioner and its approximate realization. Section~\ref{sec:experiments} illustrates
several numerical experiments which validate the theoretical findings and the practical
performance of the proposed preconditioner. Finally, contributions and future research directions are discussed
in Section~\ref{sec:conclusions}.

\section{Problem and operator form}\label{sec:model_problem}

\subsection{Notation}\label{subsec:notation}
We start by fixing some notation. Given an open and bounded domain $D$, we denote by $L^2(D)$ the
space of square integrable functions on $D$, endowed with the norm $\norm{\cdot}_{0,D}$ induced by
the inner product $\inner{\cdot}{\cdot}_{D}$. Sobolev spaces are denoted by $W^{s,p}(D)$, where
$s\in\mathbb{R}$ refers to the differentiability index and $p\in[1,\infty]$ is the integrability
exponent. When $p=2$, we use the notation $H^s(D)=W^{s,2}(D)$, with associated norm
$\norm{\cdot}_{s,D}$ and seminorm $\seminorm{\cdot}_{s,D}$.

We use normal font, e.g.\ $A$, to denote linear operators, while the matrix counterparts are
set in sans-serif font, e.g.\ $\mat{A}$. The adjoint of an
operator $A$ is denoted by $A^\ast$, and the identity operator is denoted by $I$. A calligraphic
font, e.g.\ $\Acal$, is used to denote block operators associated with saddle point
systems. A blackboard bold font, e.g.\ $\bbA$, is used to denote block matrices.

\subsection{Elliptic interface problem}
Let $\Om\subset\mathbb{R}^d$, $d\in\{2,3\}$, be an open and bounded domain with Lipschitz boundary
$\partial\Om$, and let $\Omt \subset \Om$ be a Lipschitz subdomain, compactly contained in $\Om$. We
call $\Om$ the \emph{background domain} and $\Omt$ the \emph{immersed domain}.
Setting $\Omo=\Om\setminus\overline{\Omt}$, the interface
$\Gamma\coloneqq\overline{\Omo}\cap\overline{\Omt}$ is Lipschitz continuous and satisfies
$\overline{\Gamma}\cap\partial\Om=\emptyset$. We consider the following \emph{elliptic interface problem} with a jump in the coefficients.
\begin{problem}\label{prob:elliptic_interface}
For $i=1,2$, given forcing terms $f_i\colon\Om_i\to\mathbb{R}$ and positive coefficients
$\beta_i$, find $u_i\colon\Om_i\to\mathbb{R}$ such that
\begin{equation}\label{eqn:strong_form}
    \left\{
    \begin{alignedat}{2}
        -\operatorname{div}\bigl(\beta_i\nabla u_i\bigr) & = f_i &  & \quad\text{in }\Om_i,\ i=1,2, \\
        u_1                                              & = u_2 &  & \quad\text{on }\Gamma,        \\
        \beta_1\nabla u_1\cdot\mathbf{n}_1               & = -\beta_2\nabla u_2\cdot\mathbf{n}_2
        &       &  \quad\text{on }\Gamma,         \\
        u_1                                              & = 0   &  & \quad\text{on }\partial\Om.
    \end{alignedat}
    \right.
\end{equation}
\end{problem}
In the above problem, $\mathbf{n}_i$ ($i=1,2$) denotes the unit vector normal to $\Gamma$, pointing
out of $\Om_i$. The two transmission conditions on $\Gamma$ enforce the continuity of the solution
and of the co-normal derivatives across the interface.

For the sake of simplicity, we assume $\betaone$ and $\rtwo$ to be positive constants throughout the paper,
and we assume that $\gap\coloneqq\rtwo-\betaone>0$.

\subsection{Fictitious domain reformulation}
We consider the fictitious domain reformulation of Problem~\ref{prob:elliptic_interface}
introduced in~\cite{AURICCHIO201536,BOFFI20142151}. The coefficient $\betaone$ and the forcing
term $f_1$, originally defined in $\Omo$, are extended to the whole background domain $\Om$. We
denote by $f$ the extension of $f_1$, so that $f|_{\Omo}=f_1$. Similarly, $u$ denotes the
extension of $u_1$ to $\Om$, satisfying $u|_{\Omo}=u_1$. The extended solution $u$ is then
required to match $u_2$ in the immersed domain $\Omt$, while the fictitious contribution
introduced by the extension is subtracted at the variational level on $\Omt$. To this end, we introduce the spaces
\[
    V=\Hoo,\qquad V_2=\Hot,\qquad \Lam=[\Hot]^\ast ,
\]
and the \emph{restriction} operator $\R:\Hoo\to\Hot$, defined by $\R u = u|_{\Omt}$.

Throughout, the angle bracket $\dual{\cdot}{\cdot}$ denotes the duality pairing whose spaces are understood
from its arguments: when the first argument is a multiplier in $\Lambda$, it is the $\Lambda$--$V_2$ duality on $\Omega_2$. When it pairs operators
with a test function, it is the $V$--$V^\ast$ (or $V_2$--$V_2^\ast$) duality on $\Omega$ or $\Omega_2$, as appropriate.

The matching condition $\R u=u_2$ is imposed weakly through a Lagrange multiplier $\lambda\in\Lam$, distributed
over the whole immersed domain $\Omt$, leading to the following weak formulation. We remark that other
choices for the multiplier space $\Lam$ and for the coupling term are possible~\cite{BoffiGastaldi2024}.
\begin{problem}\label{prob:weak}
Given $f\in L^2(\Omega)$ and $f_2 \in L^2(\Omega_2)$, find $(u,u_2,\lambda)\in V\times V_2\times\Lam$ such that
\begin{align}
    \betaone\inner{\nabla u}{\nabla v}_\Om + \dual{\lambda}{\R v}
                           & = \inner{f}{v}_\Om,          & \forall v   & \in V,\label{eq:aur1}   \\
    \gap\inner{\nabla u_2}{\nabla v_2}_{\Omt} - \dual{\lambda}{v_2}
                           & = \inner{f_2-f}{v_2}_{\Omt}, & \forall v_2 & \in V_2,\label{eq:aur2} \\
    \dual{\mu}{\R u - u_2} & = 0,                         & \forall\mu  & \in\Lam.\label{eq:aur3}
\end{align}
\end{problem}
Problem~\ref{prob:weak} is the standard weak formulation with fictitious domain approach of the
original interface Problem~\ref{prob:elliptic_interface} (see~\cite{AURICCHIO201536} for its derivation). The operator form of \eqref{eq:aur1}--\eqref{eq:aur3} is:
\begin{equation}\label{eq:opeq}
    \Acal
    \begin{bmatrix} u\\ u_2 \\ \lambda\end{bmatrix}
    =
    \begin{bmatrix}
        A_1 & 0   & \R^\ast \\
        0   & A_2 & -I      \\
        \R  & -I  & 0
    \end{bmatrix}
    \begin{bmatrix} u\\ u_2 \\ \lambda\end{bmatrix}
    =
    \begin{bmatrix} f\\ f_2-f \\ 0\end{bmatrix},
    \qquad
    \begin{aligned}
        A_1 & = \betaone(-\Delta^{D}_{\Om}), \\
        A_2 & = \gap(-\Delta^{N}_{\Omt}),
    \end{aligned}
\end{equation}where $\Delta^{D}_{\Om}$ and $\Delta^{N}_{\Omt}$ denote the Laplace operators with homogeneous Dirichlet and Neumann boundary conditions, respectively. Defining $\W \coloneqq V\times V_2$, $B\coloneqq [\,\R\ \ -I\,]:\W\to V_2$, and $A \coloneqq \operatorname{diag}(A_1,A_2)$, $\Acal$ has a standard saddle point structure.
The associated bilinear and linear forms are then
\begin{align}
    a((u,u_2),(v,v_2)) & = \betaone \inner{\nabla u}{\nabla v}_\Om + \gap\inner{\nabla u_2}{\nabla v_2}_{\Omt}, \\
    b(\lambda,(v,v_2)) & = \dual{\lambda}{\R v - v_2},                                                          \\
    F((v,v_2))         & = \inner{f}{v}_\Om + \inner{f_2-f}{v_2}_{\Omt},
\end{align}

Well-posedness for this saddle point problem follows from the classical saddle point theory~\cite{boffi2013mixed}, upon verifying the
ellipticity of $a(\cdot,\cdot)$ on the kernel of $b(\cdot,\cdot)$ and the inf-sup condition. Detailed proofs can be
found in~\cite{AURICCHIO201536,BOFFI20142151,BoffiGastaldi2024}. Details about the regularity of the solution can be found in~\cite{nicaise,AURICCHIO201536}.

We now introduce norms on the spaces $\W$ and $\Lam$ that will be used to build a
parameter-robust preconditioner.

\section{Operator based preconditioning}\label{sec:operator_prec}

We choose norms on $W$ and $\Lam$ that will make the Brezzi conditions of $\Acal$ independent
of $\betaone,\rtwo$, which we use to construct a preconditioner using the corresponding Riesz maps. We start by defining the following
coefficient-dependent norm on the immersed domain $\Omega_2$:
\begin{equation}\label{eq:imm-norm}
    \norm{q}_{\ast}^2 \coloneqq \gap\,\seminorm{q}_{1,\Omt}^2 + \betaone\,\norm{q}_{0,\Omt}^2 ,
    \qquad q\in\Hot,
\end{equation}
and the unweighted one $\norm{q}_{V_2}^2\coloneqq\seminorm{q}_{1,\Omt}^2+\norm{q}_{0,\Omt}^2$. Its dual norm, denoted by $\norm{\mu}_{V_2^{\ast}}$, is naturally
defined as
\begin{equation}\label{eqn:flat}
    \norm{\mu}_{V_2^{\ast}}\coloneqq\sup_{0\ne q \in V_2} \frac{\dual{\mu}{q}}{\norm{q}_{V_2}}.
\end{equation}On the space $W = V \times V_2$, we define the following norm
\begin{equation}\label{eq:Wnorm}
    \norm{(u,u_2)}_\W^2
    \coloneqq
    \betaone\,\seminorm{u}_{1,\Om}^2
    +
    \gap\,\seminorm{u_2}_{1,\Omt}^2 + \betaone\,\norm{u_2}_{0,\Omt}^2 = \betaone\,\seminorm{u}_{1,\Om}^2 + \norm{u_2}_\ast^2,
\end{equation}while we define the multiplier norm on $\Lam$ as the sum of background and immersed contributions:
\begin{equation}\label{eq:Qnorm}
    \begin{aligned}
        \norm{\lambda}_\Q^2
         & \coloneqq
        \tfrac{1}{\betaone}\,\dual{\lambda}{\R(-\Delta^{D}_\Om)^{-1}\R^\ast\lambda}%
        +
        \dual{\lambda}{(\gap(-\Delta^{N}_{\Omt})+\betaone I)^{-1}\lambda}%
        \\
         & = \norm{\lambda}_{\Q,\mathrm{backgr}}^2 + \norm{\lambda}_{\Q,\mathrm{imm}}^2.
    \end{aligned}
\end{equation}

We notice that the natural candidates for the first two diagonal blocks of the preconditioner are
$\bigl(\betaone(-\Delta^D_\Om)\bigr)^{-1}$ and $\bigl(\gap(-\Delta^N_{\Omt})+\betaone I\bigr)^{-1}$, which
are the Riesz maps associated with the weighted norms on $V$ and $V_2$, respectively.

Clearly, both operators will lead to symmetric and positive definite (SPD) matrices upon discretization. In particular, the shift
$\betaone I$ removes the singularity from $A_2$, which is positive semidefinite\footnote{Indeed, no boundary conditions are imposed in $V_2$.}. The choice of the multiplier block is crucial in order to obtain a parameter-robust preconditioner. We present some useful results which will be needed for the analysis.

\subsection{Preliminaries}
We start by recalling some fundamental results from~\cite{evans2022partial}, applied to the present setting.

\begin{lemma}[Poincar\'e]\label{lem:poincare}
    There is $C_P=C_P(\Om)$ with $\norm{v}_{0,\Om}\le C_P\seminorm{v}_{1,\Om}$ for all
    $v \in \Hoo$.
    In particular, for $u\in\Hoo$,
    \[
        \seminorm{\R u}_{1,\Omt}\le \seminorm{u}_{1,\Om},
        \qquad
        \norm{\R u}_{0,\Omt}\le\norm{u}_{0,\Om}\le C_P\seminorm{u}_{1,\Om},
    \]
    so that $\norm{\R u}_{V_2}^2\le (1+C_P^2)\seminorm{u}_{1,\Om}^2$.
\end{lemma}

\begin{lemma}[Extension theorem]\label{lem:ext}
    There is a bounded linear extension $\E:\Hot\to\Hoo$ with $\R\E=\mathrm{Id}$ and such that
    \[
        \seminorm{\E q}_{1,\Om}\le C_E\,\norm{q}_{V_2},\qquad\forall q\in\Hot,
    \]
    where $C_E=C_E(\Om,\Omt) \geq1$ depends only on the geometry. Such an extension $\E$ exists because $\Omt\Subset\Om$.
\end{lemma}
To see that $C_E \geq 1$, consider the oscillatory sequence $\{ q_n(\boldsymbol{x}) \}_n = \frac{1}{n} \sin(n x_1)$ in $\Omt$, normalized
so that $|q_n|_{1, \Omt}=1$, while $\|q_n\|_{0, \Omt} \to_{n\to\infty} 0$. Then $\norm{q_n}_{V_2} \to_{n\to\infty} 1$ and, since $Eq_n = q_n$
on $\Omega_2$, $\seminorm{\E q_n}_{1,\Om} \geq \seminorm{q_n}_{1,\Omt} \!=\! 1$, showing that $C_E \geq (1+ \|q_n\|_{0,\Omt}^2)^{-1/2} \to_{n\to\infty} 1$.

\begin{lemma}\label{lem:schur}
    For every $\mu\in\Lam$, it holds
    \[
        C_E^{-2}\,\norm{\mu}_{V_2^{\ast}}^2
        \le
        \sup_{0\ne u\in\Hoo}\frac{\dual{\mu}{\R u}^2}{\seminorm{u}_{1,\Om}^2}
        \le
        (1+C_P^2)\,\norm{\mu}_{V_2^{\ast}}^2.
    \]
\end{lemma}

\begin{proof}
    Let $u\in\Hoo$. Lemma~\ref{lem:poincare} gives $\norm{\R u}_{V_2}^2\le(1+C_P^2)\seminorm{u}_{1,\Om}^2$, so
    \[\dual{\mu}{\R u}^2\le\norm{\mu}_{V_2^{\ast}}^2\norm{\R u}_{V_2}^2 \le(1+C_P^2)\norm{\mu}_{V_2^{\ast}}^2\seminorm{u}_{1,\Om}^2,\]where the first inequality exploits the definition of dual norm. Conversely, for $0 \ne q\in\Hot$, we set $u=\E q$. Then, $\R u=q$ and, by
    Lemma~\ref{lem:ext}, $\seminorm{u}_{1,\Om}\le C_E\norm{q}_{V_2}$, whence
    $\dual{\mu}{q}^2/\norm{q}_{V_2}^2\le C_E^2\,\dual{\mu}{\R u}^2/\seminorm{u}_{1,\Om}^2$. Taking suprema gives the claim.
\end{proof}

\begin{lemma}\label{lem:dualnorm}
    Let $(q,p)_\ast\coloneqq\gap\inner{\nabla q}{\nabla p}_{\Omt}+\betaone\inner{q}{p}_{\Omt}$ induce
    $\norm{\cdot}_\ast$, and let $S\coloneqq\gap(-\Delta^N_{\Omt})+\betaone I:V_2\to V_2^\ast$, which is
    an isomorphism since $\betaone>0$. For $0 \ne \mu\in\Lam$, let $r \coloneqq S^{-1}\mu\in V_2$ be its Riesz
    representative, i.e.\ $\dual{\mu}{p}=(r,p)_\ast$ for all $p\in V_2$. Then
    \[
        \norm{r}_\ast=\norm{\mu}_{\Q,\mathrm{imm}}
        =\sup_{0\ne\eta\in V_2}\frac{\dual{\mu}{\eta}}{\norm{\eta}_\ast}.
    \]
\end{lemma}

\begin{proof}
    First, $\norm{r}_\ast^2=(r,r)_\ast=\dual{\mu}{r}=\dual{\mu}{S^{-1}\mu}
        =\norm{\mu}_{\Q,\mathrm{imm}}^2$ by \eqref{eq:Qnorm}. For the supremum: for every
    $\eta \in V_2$, $\eta\ne0$, the Riesz identity and Cauchy--Schwarz give
    $\dual{\mu}{\eta}=(r,\eta)_\ast\le\norm{r}_\ast\norm{\eta}_\ast$, hence
    $\sup_{0\ne\eta\in V_2} \dual{\mu}{\eta}/\norm{\eta}_\ast\le\norm{r}_\ast$.

    Testing at $\eta=r$ gives $\dual{\mu}{r}/\norm{r}_\ast=(r,r)_\ast/\norm{r}_\ast=\norm{r}_\ast$. Hence the supremum is attained and equals $\norm{r}_\ast$.
\end{proof}

\subsection{Robust preconditioning for the elliptic interface problem}

We can now state the well-posedness of the problem on $W \times \Lam$ using the parameter-dependent norms previously defined.
\begin{theorem}\label{thm:cont}
    Let $\gap>0$ and equip $W,$ $\Lam$ with the norms~\eqref{eq:Wnorm} and~\eqref{eq:Qnorm}, respectively. Then
    $\Acal \colon W\times \Lam \to(W\times \Lam)^\ast$ is an isomorphism, and its condition number is bounded independently of $\betaone,\rtwo$ and of their contrast.
\end{theorem}

\begin{proof}
    Following~\cite[Appendix A]{prec1d2d}, we verify the four Brezzi conditions~\cite{BrezziFortin}: boundedness of
    $a(\cdot,\cdot)$ and $b(\cdot,\cdot)$, coercivity of $a(\cdot,\cdot)$ on $\ker b$, and the inf-sup condition for $b(\cdot,\cdot)$.

    First, we establish the boundedness of $a(\cdot,\cdot)$. Let $w=(u,u_2) \in W$.
    Since $A$ is block-diagonal, $\dual{Aw}{w}=\betaone\seminorm{u}_{1,\Om}^2+\gap\seminorm{u_2}_{1,\Omt}^2\le\norm{w}_W^2$. By Cauchy--Schwarz, $\dual{Aw}{w'}\le\norm{w}_W\norm{w'}_W$.

    Next, we prove coercivity on $\ker b$. For $w\in\ker b$, we have $u_2=\R u$, so by Lemma~\ref{lem:poincare}:
    $\seminorm{u_2}_{1,\Omt}\le\seminorm{u}_{1,\Om}$ and
    $\norm{u_2}_{0,\Omt}\le C_P\seminorm{u}_{1,\Om}$. Then,
    $\dual{Aw}{w}=\betaone\seminorm{u}_{1,\Om}^2+\gap\seminorm{u_2}_{1,\Omt}^2,$
    while
    \[
        \norm{w}_W^2
        = \betaone\seminorm{u}_{1,\Om}^2+\gap\seminorm{u_2}_{1,\Omt}^2
        + \betaone\norm{u_2}_{0,\Omt}^2
        \le \betaone\seminorm{u}_{1,\Om}^2+\gap\seminorm{u_2}_{1,\Omt}^2
        + \betaone C_P^2\seminorm{u}_{1,\Om}^2 .
    \]
    Since $\betaone \norm{u}_{0,\Omega_2}^2 \leq  \betaone C_P^2 \seminorm{u}_{1,\Omega}^2 \leq C_P^2 \dual{Aw}{w}$, it follows
    \[
        \norm{w}_W^2
        = \dual{Aw}{w} + \betaone \norm{u_2}_{0,\Omega_2}^2 \leq \dual{Aw}{w} +  C_P^2 \dual{Aw}{w} = (1+C_P^2) \dual{Aw}{w},
    \]which gives $\dual{Aw}{w} \geq \frac{1}{1+C_P^2} \norm{w}_{W}^2.$

    To establish the boundedness of $b(\cdot,\cdot)$, we recall the two terms in the multiplier norm~\eqref{eq:Qnorm},
    and that $\norm{\mu}_\Q^2=\norm{\mu}_{\Q,\mathrm{backgr}}^2+\norm{\mu}_{\Q,\mathrm{imm}}^2$.

    We start by bounding the term $\norm{\mu}_{\Q,\mathrm{imm}}$. Let $(q,p)_\ast\coloneqq\gap\inner{\nabla q}{\nabla p}_{\Omt}+\betaone\inner{q}{p}_{\Omt}$ be the inner
    product inducing the norm $\norm{\cdot}_\ast$ in~\eqref{eq:imm-norm}, and let
    $S\coloneqq\gap(-\Delta^N_{\Omt})+\betaone I:V_2\to V_2^\ast$. Let $r$ be the Riesz representative of $\mu$, then for all $p \in V_2$:
    \[
        \dual{\mu}{p}=(r,p)_\ast = \alpha (\nabla r, \nabla p)_{\Omega_2} + \betaone(r,p)_{\Omega_2} = \dual{Sr}{p},
    \]from which it follows $r = S^{-1} \mu$. Consequently, the following identity between the norms of $r$ and $\mu$ holds:
    $\norm{r}_{\ast}^2=(r,r)_{\ast}=\dual{\mu}{r} = \dual{\mu}{S^{-1}\mu} =
        \norm{\mu}_{\Q,\mathrm{imm}}^2.$
    Therefore, using the Riesz identity with $p=u_2$ and then Cauchy--Schwarz inequality, we have:
    \begin{equation}\label{eq:Bc-imm}
        |\dual{\mu}{u_2}| = |(r,u_2)_{\ast}| \le \norm{r}_{\ast} \norm{u_2}_{\ast}
        \le\ \norm{\mu}_{\Q,\mathrm{imm}}\,\norm{u_2}_\ast.
    \end{equation}

    To bound $\norm{\mu}_{\Q,\mathrm{backgr}}$, we first notice that the functional $v\mapsto\dual{\mu}{\R v}$ is bounded on $\Hoo$ endowed with the inner product $\betaone\inner{\nabla\cdot}{\nabla\cdot}_\Om$. Indeed, applying the definition of dual norm and Lemma~\ref{lem:poincare}:
    \[\lvert\dual{\mu}{\R v}\rvert\le\norm{\mu}_{V_2^\ast}\norm{\R v}_{V_2}
        \le\norm{\mu}_{V_2^\ast}\sqrt{1+C_P^2}\,\seminorm{v}_{1,\Om}.\]
    Hence there is a unique Riesz representative $z=z(\mu)\in\Hoo$ which satisfies the problem
    \begin{equation}\label{eq:Bc-star}
        \betaone\inner{\nabla z}{\nabla v}_\Om=\dual{\mu}{\R v}\qquad\forall v\in\Hoo,
    \end{equation}
    or $z=\tfrac1{\betaone}(-\Delta^D_\Om)^{-1}\R^\ast\mu$, so that
    $(-\Delta^D_\Om)^{-1}\R^\ast\mu=\betaone z$ and
    $\R(-\Delta^D_\Om)^{-1}\R^\ast\mu=\betaone\,\R z$. Substituting into the definition of $\norm{\cdot}_{\Q,\mathrm{backgr}}$, and testing~\eqref{eq:Bc-star} with $v=z$, we
    obtain:
    \begin{equation}\label{eq:Bc-extnorm}
        \begin{aligned}
            \norm{\mu}_{\Q,\mathrm{backgr}}^2
             & =\tfrac1{\betaone}\dual{\mu}{\R(-\Delta^D_\Om)^{-1}\R^\ast\mu}
            =\tfrac1{\betaone}\dual{\mu}{\betaone\,\R z}                      \\
             & =\dual{\mu}{\R z}
            \overset{\eqref{eq:Bc-star}}{=}\betaone\,\seminorm{z}_{1,\Om}^2 .
        \end{aligned}
    \end{equation}
    Take now an arbitrary $u\in\Hoo$: applying \eqref{eq:Bc-star} with $v=u$, the
    Cauchy--Schwarz inequality for $\inner{\cdot}{\cdot}_\Om$, and~\eqref{eq:Bc-extnorm}, yields a bound for $\dual{\mu}{\R u}$:
    \begin{equation}\label{eq:Bc-ext}
        |\dual{\mu}{\R u}|
        = |\betaone\inner{\nabla z}{\nabla u}_\Om|
        \le\betaone\,\seminorm{z}_{1,\Om}\seminorm{u}_{1,\Om}
        =\norm{\mu}_{\Q,\mathrm{backgr}}\,\bigl(\sqrt{\betaone}\,\seminorm{u}_{1,\Om}\bigr).
    \end{equation}

    We can combine the previous two bounds to show that $B$ is bounded. Given $w=(u,u_2)\in\W$ and $\mu\in\Q$, using $\dual{Bw}{\mu}=\dual{\mu}{\R u}-\dual{\mu}{u_2}$
    together with bounds~\eqref{eq:Bc-imm} and~\eqref{eq:Bc-ext},
    \[
        \dual{Bw}{\mu}
        \le \bigl\lvert\dual{\mu}{\R u}\bigr\rvert + \bigl\lvert\dual{\mu}{u_2}\bigr\rvert
        \le \norm{\mu}_{\Q,\mathrm{backgr}}\,\bigl(\sqrt{\betaone}\,\seminorm{u}_{1,\Om}\bigr)
        +\norm{\mu}_{\Q,\mathrm{imm}}\,\norm{u_2}_\ast.
    \]
    Applying the Cauchy--Schwarz inequality in $\mathbb{R}^2$ to the vectors
    $\bigl(\norm{\mu}_{\Q,\mathrm{backgr}},\norm{\mu}_{\Q,\mathrm{imm}}\bigr)$ and
    $\bigl(\sqrt{\betaone}\,\seminorm{u}_{1,\Om},\norm{u_2}_\ast\bigr)$, we get the desired bound
    \[
        \dual{Bw}{\mu}
        \le\Bigl(\norm{\mu}_{\Q,\mathrm{backgr}}^2+\norm{\mu}_{\Q,\mathrm{imm}}^2\Bigr)^{1/2}
        \Bigl(\betaone\seminorm{u}_{1,\Om}^2+\norm{u_2}_\ast^2\Bigr)^{1/2}
        =\norm{\mu}_\Q\,\norm{w}_\W ,
    \]
    where the final equality uses $\norm{\mu}_\Q^2=\norm{\mu}_{\Q,\mathrm{backgr}}^2
        +\norm{\mu}_{\Q,\mathrm{imm}}^2$ and the identity
    $\norm{w}_\W^2=\betaone\seminorm{u}_{1,\Om}^2+\norm{u_2}_\ast^2$ from \eqref{eq:Wnorm}.

    Finally, we verify the inf-sup condition. Assume $0 \ne \mu$, and let $z\in\Hoo$ and $r\in V_2$ be the background and immersed Riesz
    representatives of $\mu$ from \eqref{eq:Bc-star} and Lemma~\ref{lem:dualnorm}, and take
    the test pair $w=(z,-r)\in\W$. Then $\R z-(-r)=\R z+r$, so by \eqref{eq:Bc-extnorm} and
    Lemma~\ref{lem:dualnorm},
    \[
        \dual{Bw}{\mu}=\dual{\mu}{\R z}+\dual{\mu}{r}
        =\norm{\mu}_{\Q,\mathrm{backgr}}^2+\norm{\mu}_{\Q,\mathrm{imm}}^2=\norm{\mu}_\Q^2 ,
    \]
    \[
        \norm{w}_\W^2=\betaone\seminorm{z}_{1,\Om}^2+\norm{r}_\ast^2
        =\norm{\mu}_{\Q,\mathrm{backgr}}^2+\norm{\mu}_{\Q,\mathrm{imm}}^2=\norm{\mu}_\Q^2 .
    \]
    Hence
    \[
        \sup_{0 \ne w'\in\W}\frac{\dual{Bw'}{\mu}}{\norm{w'}_\W}
        \ \ge\ \frac{\dual{Bw}{\mu}}{\norm{w}_\W}
        =\frac{\norm{\mu}_\Q^2}{\norm{\mu}_\Q}=\norm{\mu}_\Q.
    \]Dividing by $\norm{\mu}_\Q$ and taking the infimum over $0 \ne \mu\in\Q$ gives the thesis.
\end{proof}

\begin{corollary}\label{cor:Aext}
    Let $\norm{\cdot}_{V_2^{\ast}}$ be defined in~\eqref{eqn:flat}. Then, the following equivalence holds for any $0\ne \mu\in\Q$:
    \[
        (1+C_P^2)^{-1/2}\,\norm{\mu}_{\Q,\mathrm{backgr}}
        \ \le\
        \tfrac{1}{\sqrt{\betaone}}\,\norm{\mu}_{V_2^{\ast}}
        \ \le\
        C_E\,\norm{\mu}_{\Q,\mathrm{backgr}} .
    \]
\end{corollary}

\begin{proof}
    Let $z=z(\mu)\in\Hoo$ solve $\betaone\inner{\nabla z}{\nabla v}_\Om=\dual{\mu}{\R v}$ for all
    $v\in\Hoo$. Then, $\dual{\mu}{\R u}=\betaone\inner{\nabla z}{\nabla u}_\Om$ and, by
    Cauchy--Schwarz inequality we have

    $$\sup_{0 \ne u \in H_0^1(\Omega)}\frac{\dual{\mu}{\R u}^2}{\seminorm{u}_{1,\Om}^2}\le\sup_{0 \ne u \in H_0^1(\Omega)}\frac{\betaone^2 \seminorm{z}_{1,\Om}^2 \seminorm{u}_{1,\Om}^2}{\seminorm{u}_{1,\Om}^2}=\betaone^2\seminorm{z}_{1,\Om}^2.$$

    We now observe that the supremum is attained at $u=z$, hence
    $$\sup_{0 \ne u \in H_0^1(\Omega)}\dual{\mu}{\R u}^2/\seminorm{u}_{1,\Om}^2=\betaone^2\seminorm{z}_{1,\Om}^2
        =\betaone\norm{\mu}_{\Q,\mathrm{backgr}}^2,$$ using
    $\norm{\mu}_{\Q,\mathrm{backgr}}^2=\betaone\seminorm{z}_{1,\Om}^2$ from
    \eqref{eq:Bc-extnorm}. Substituting this into Lemma~\ref{lem:schur} gives
    $$C_E^{-2}\norm{\mu}_{V_2^{\ast}}^2\le\betaone\norm{\mu}_{\Q,\mathrm{backgr}}^2\le
        (1+C_P^2)\norm{\mu}_{V_2^{\ast}}^2.$$ Finally, dividing by $\betaone>0$ and rearranging each side
    yields the two inequalities.
\end{proof}

\subsection{The multiplier block}\label{sec:multiplier_block}

In this section we show that the multiplier block $N_\Q$, which is  by definition the sum of the two
operators in~\eqref{eq:Qnorm},
\[
    N_\Q\;\coloneqq\;\tfrac1{\betaone}\R(-\Delta^D_\Om)^{-1}\R^\ast
    \;+\;\bigl(\gap(-\Delta^N_{\Omt})+\betaone I\bigr)^{-1},
    \qquad \norm{\mu}_\Q^2=\dual{\mu}{N_\Q\mu},
\]is spectrally equivalent, with constants independent of the contrast between subdomains, to a \emph{single} operator which depends on $-\Delta^N_{\Omt}$ and $I$ only. To this
end, we replace the first summand in $N_\Q$ with a spectrally equivalent operator involving $-\Delta^N_{\Omt}$.

\begin{definition}
    Let $\{(\eta_k,\phi_k)\}_{k\ge1}$ be the $L^2(\Omt)$-orthonormal eigenbasis of the
    Neumann Laplacian,
    \[
        (-\Delta^N_{\Omt})\phi_k=\eta_k\,\phi_k,\qquad 0\le\eta_k\to\infty,\qquad
        \inner{\phi_k}{\phi_j}_{\Omt}=\delta_{kj}.
    \]
    In the following calculations, $q$ denotes an arbitrary linear combination $q= \sum_{k=1}^m c_k \phi_k$ of the Neumann
    eigenfunctions, with $c_k = \inner{q}{\phi_k}_{\Omt}$. Then, $\norm{q}_{0,\Omt}^2=\sum_{k=1}^m c_k^2$ and
    $\seminorm{q}_{1,\Omt}^2=\sum_{k=1}^m \eta_k c_k^2$.
\end{definition}
Every operator built from $-\Delta^N_{\Omt}$ and $I$ in this section is diagonal in this
basis. Indeed, an operator $F=f(-\Delta^N_{\Omt})$, defined by $F\phi_k=f(\eta_k)\phi_k$ for a
suitable scalar function $f$, has $k$-th eigenvalue $f(\eta_k)$ and
$\dual{q}{Fq}=\sum_k f(\eta_k)c_k^2$. To give an example, the operator $\gap(-\Delta^N_{\Omt})+\betaone I$ has
eigenvalues $\gap\eta_k+\betaone$, and its inverse has eigenvalues $(\gap\eta_k+\betaone)^{-1}$.

We say two operators $F$ and $G$ are \emph{spectrally equivalent}, and write $F\equeq G$, when
\[
    c\dual{q}{Gq}\le\dual{q}{Fq}\le C\dual{q}{Gq} \quad\text{for all }q,
\]
for some $c,C>0$ which do not depend on discretization parameters. It is immediate to check that $\equeq$ is
an equivalence relation. The following Lemma provides a convenient criterion to check spectral equivalence.

\begin{lemma}\label{lem:symcalc}
    Let $F=f(-\Delta^N_{\Omt})$ and $G=g(-\Delta^N_{\Omt})$ be diagonalized by the common
    eigenbasis $\{\phi_k\}_k$ with eigenvalues $f(\eta_k),g(\eta_k)>0$, and let $c,C>0$ be positive
    constants. Then,
    \[
        c\,\dual{q}{Gq}\le \dual{q}{Fq}\le C\,\dual{q}{Gq}\quad\text{for all }q
        \qquad\Longleftrightarrow\qquad
        c\le\frac{f(\eta_k)}{g(\eta_k)}\le C\quad\text{for all }k .
    \]
    That is, $F\equeq G$ with constants $c,C$ iff the ratio of eigenvalues of $F$
    and $G$ lies in $[c,C]$.
\end{lemma}

\begin{proof}
    Writing $q=\sum_k c_k\phi_k$, one has $\dual{q}{Fq}=\sum_k f(\eta_k)c_k^2$ and
    $\dual{q}{Gq}=\sum_k g(\eta_k)c_k^2$.

    From $c\le f(\eta_k)/g(\eta_k)\le C$ for all $k$, we multiply by $g(\eta_k)c_k^2$ and sum over $k$ to obtain
    \[
        c\dual{q}{Gq}=\sum_k c\,g(\eta_k)c_k^2\le\sum_k f(\eta_k)c_k^2=\dual{q}{Fq}\le\sum_k C\,g(\eta_k)c_k^2=C\dual{q}{Gq}.
    \]

    Conversely, testing on $q=\phi_k$ and using the orthonormality of the basis, it follows $c\le f(\eta_k)/g(\eta_k)\le C$.
\end{proof}

Thus, for operators diagonalized by the basis $\{\phi_k\}_k$, spectral equivalence reduces
to a two-sided bound on the ratio of eigenvalues, with the same constants. We now show that the background term, $\tfrac1{\betaone}\R(-\Delta^D_\Om)^{-1}\R^\ast$, is spectrally equivalent
to a function of $-\Delta^N_{\Omt}$ and $I$, with constants depending only on the geometries $\Om$ and $\Omt$.
\begin{lemma}\label{lem:extreduce}
    As operators on $\Lam$,
    \[
        \tfrac1{\betaone}\,\R(-\Delta^D_\Om)^{-1}\R^\ast
        \ \equeq\
        \tfrac1{\betaone}\bigl(-\Delta^N_{\Omt}+I\bigr)^{-1},
    \]
    with equivalence constants $C_E^{-2}$ and $1+C_P^2$, depending only on $\Om,\Omt$.
\end{lemma}

\begin{proof}
    The unweighted norm is $\norm{q}_{V_2}^2 \coloneqq \seminorm{q}_{1,\Omt}^2+\norm{q}_{0,\Omt}^2
        =\dual{(-\Delta^N_{\Omt}+I)q}{q}$. Then, repeating the proof of
    Lemma~\ref{lem:dualnorm} with $(q,p)_{V_2} \coloneqq \inner{\nabla q}{\nabla p}_{\Omt}+\inner{q}{p}_{\Omt}$, we
    obtain
    \begin{equation}\label{eq:flatdual}
        \norm{\mu}_{V_2^{\ast}}^2=\dual{\mu}{(-\Delta^N_{\Omt}+I)^{-1}\mu}.
    \end{equation}

    Moreover, by the computation in the proof of Corollary~\ref{cor:Aext},
    \begin{equation}\label{eq:extdual}
        \sup_{0\ne u\in\Hoo}\frac{\dual{\mu}{\R u}^2}{\seminorm{u}_{1,\Om}^2}
        =\dual{\mu}{\R(-\Delta^D_\Om)^{-1}\R^\ast\mu}.
    \end{equation}
    Substituting the identities~\eqref{eq:flatdual} and \eqref{eq:extdual} into Lemma~\ref{lem:schur} gives, for all $\mu\in\Lam$,
    \[
        C_E^{-2}\,\dual{\mu}{(-\Delta^N_{\Omt}+I)^{-1}\mu}
        \ \le\
        \dual{\mu}{\R(-\Delta^D_\Om)^{-1}\R^\ast\mu}
        \ \le\
        (1+C_P^2)\,\dual{\mu}{(-\Delta^N_{\Omt}+I)^{-1}\mu},
    \]
    which gives the claim after dividing all terms by $\betaone>0$.
\end{proof}

Thanks to Lemma~\ref{lem:extreduce}, we replace the background term in $N_\Q$ with its spectrally equivalent shifted Neumann Laplacian inverse, which gives
the operator
\begin{equation}\label{eq:G}
    G\coloneqq\tfrac1{\betaone}(-\Delta^N_{\Omt}+I)^{-1}
    +\bigl(\gap(-\Delta^N_{\Omt})+\betaone I\bigr)^{-1}.
\end{equation}In particular, it is immediate to check that
\begin{equation}\label{eq:NG}
    N_\Q\ \equeq\ G,
\end{equation}with the same constants $C_E^{-2},\,1+C_P^2$. The operator $G$ is again diagonal in
the basis $\{\phi_k\}_k$, with $k$-th eigenvalue $s(\eta_k)$, where
\begin{equation}\label{eq:symbol}
    s(\eta)\coloneqq\frac{1}{\betaone(\eta+1)}\;+\;\frac{1}{\gap\eta+\betaone},
    \qquad \eta\ge0.
\end{equation}

We can finally state the spectral equivalence of the multiplier block $N_\Q$ with
a single operator $N$, uniformly in the contrast between coefficients. To do so, we exploit the condition
in Lemma~\ref{lem:symcalc}.

\begin{theorem}\label{thm:multiplier_block}
    Let $\theta\coloneqq\betaone(\rtwo-\betaone)/\rtwo$, and $N\coloneqq\theta(-\Delta^N_{\Omt})+\betaone I$. Then $N_\Q \equeq N^{-1}$,
    i.e., the multiplier block $N_\Q$ is spectrally equivalent
    to $N^{-1}=\bigl(\theta(-\Delta^N_{\Omt})+\betaone I\bigr)^{-1}$, uniformly in the jump,
    with equivalence constants depending only on $\Om$, $\Omt$, and on the
    absolute factors of Lemma~\ref{lem:s_bound}, but not on $\betaone,\rtwo$.
\end{theorem}

\begin{proof}
    By the equivalence result in~\eqref{eq:NG}, $N_\Q\equeq G$ with constants $C_E^{-2},1+C_P^2$.
    The operators $G$ and $N^{-1}$ are both diagonal in the basis $\{\phi_k\}_k$, with
    $k$-th eigenvalues $s(\eta_k)$ and $(\theta\eta_k+\betaone)^{-1}$, respectively.

    Lemma~\ref{lem:s_bound} states the following bound for $s(\eta)^{-1}$:
    \[
        \tfrac12 (\theta\eta+\betaone)\le s(\eta)^{-1}\le \theta\eta+\betaone.
    \]Taking now reciprocals, we obtain
    \begin{equation}\label{eq:ratio}
        1 \le\ \frac{s(\eta_k)}{(\theta\eta_k+\betaone)^{-1}}\ \le\ 2
        \qquad\text{for all }k .
    \end{equation}
    By the condition in Lemma~\ref{lem:symcalc}, the bound in~\eqref{eq:ratio} implies
    $G \equeq N^{-1}$ with constants $1$ and $2$. Thesis follows from the transitivity of $\equeq$.
\end{proof}

By Theorems~\ref{thm:cont} and~\ref{thm:multiplier_block}, and the operator preconditioning framework~\cite{MW}, a spectrally equivalent
approximation to the Riesz map preconditioner associated with the product norm on $W \times \Lam$ defined by~\eqref{eq:Wnorm}--\eqref{eq:Qnorm} is the operator
\begin{equation}\label{eq:precond}
    \Bcal \coloneqq \diag\left(
    \bigl(\betaone(-\Delta^D_\Om)\bigr)^{-1},
    \bigl(\gap(-\Delta^N_{\Omt})+\betaone I\bigr)^{-1},
    N %
    \right).
\end{equation}
We notice that \emph{all} the blocks of the preconditioner can be realized, at the discrete level, in an efficient and order-optimal manner through multilevel schemes. Remarkably, and in contrast with the codimension-one framework
established in~\cite{prec1d2d}, no fractional operators are involved in this case for the multiplier term in the last block.
\begin{remark}[Comparison with codimension-one solvers]
    In previous collaborative work~\cite{BENZI2026118522}, the author proposed an augmented Lagrangian preconditioner for the Poisson and Stokes fictitious domain
    problems, when the constraint is enforced on a codimension-one domain. In such a case, the operator preconditioning approach entails the use
    of the fractional Laplacian~\cite{prec1d2d}, which is expensive to realize as it requires the solution of several linear
    systems~\cite{haznics}. Indeed, the performance of the AL preconditioner was superior both in terms of wall-clock
    times and iteration counts.

    In contrast, the present work shows that for codimension-zero problems, the diagonal blocks of
    the preconditioner are realized by simpler differential operators, which are much easier to invert
    than the ones in the AL approach. This will be crucial for the practical performance of the preconditioner. In particular, efficient matrix-free operator evaluations, combined with
    the optimal convergence properties of geometric multigrid methods, offer an excellent avenue to achieve scalable simulations.
    These observations are backed by numerical evidence in Section~\ref{sec:experiments}.
\end{remark}

\section{Finite element discretization and discrete preconditioner}\label{sec:discretization}

We adopt the finite element discretization considered in~\cite{AURICCHIO201536,BOFFI20142151}.

We employ two \emph{independent}, shape regular mesh families for $\Omega$ and $\Omega_2$, which we denote by $\mathcal{T}_h$ and $\mathcal{T}_{2,h}$ respectively.
In the fictitious domain setting, the background domain is taken as a structured quadrilateral (or hexahedral) discretization of a $d$-dimensional box. Hence, the background domain $\Omega$ will be assumed
to be convex. In addition, we require the immersed mesh $\mathcal{T}_{2,h}$ to be quasi-uniform. We denote by $h$ and $h_2$ the mesh sizes for $\mathcal{T}_h$ and $\mathcal{T}_{2,h}$, respectively. Moreover, for a generic element $K \in \mathcal{T}_h$ or $\mathcal{T}_{2,h}$, we define $\mathcal{Q}^p(K)$, $p\geq1$, to be the space of polynomials defined on $K$ of degree at most $p$ in each variable.

After introducing finite dimensional subspaces $V_h\subset V$, $V_{2,h}\subset V_2$, and $\Lambda_h \subset \Lambda$, the discrete problem associated with the continuous Problem~\ref{prob:weak} reads as follows.

\begin{problem}
\label{prob:discrete}
Given $f \in L^2(\Omega)$ and $f_2 \in L^2(\Omega_2)$, find $(u_h,u_{2,h},\lambda_h)\in V_h\times V_{2,h}\times \Lambda_h$ such that
\begin{equation}\label{eq:discrete}
    \begin{aligned}
        \betaone\inner{\nabla u_h}{\nabla v_h}_\Om+\dual{\lambda_h}{\R v_h}         & =\inner{f}{v_h}_\Om
                                                                                    &                                & \forall v_h\in V_h,         \\
        \gap\inner{\nabla u_{2,h}}{\nabla v_{2,h}}_{\Omt}-\dual{\lambda_h}{v_{2,h}} & =\inner{f_2-f}{v_{2,h}}_{\Omt}
                                                                                    &                                & \forall v_{2,h}\in V_{2,h}, \\
        \dual{\mu_h}{\R u_h-u_{2,h}}                                                & =0
                                                                                    &                                & \forall \mu_h\in \Lambda_h.
    \end{aligned}
\end{equation}
\end{problem}

We observe that the duality pairing $\dual{\mu}{v_2}$ can be identified with the scalar product in $L^2(\Omega_2)$ if $\mu \in L^2(\Omega_2)$.
To solve Problem~\ref{prob:discrete}, we consider conforming Lagrangian elements for all the variables:
\begin{align*}
    V_h       & = \{v_h \in H_0^1(\Omega): v_h|_K \in \mathcal{Q}^1(K) \text{ for all } K \in \mathcal{T}_h\} \subset H_0^1(\Omega),             \\
    V_{2,h}   & = \{v_{2,h} \in H^1(\Omega_2): v_{2,h}|_K \in \mathcal{Q}^1(K) \text{ for all } K \in \mathcal{T}_{2,h}\} \subset H^1(\Omega_2), \\
    \Lambda_h & = V_{2,h}.
\end{align*}

The interested reader is referred to~\cite{AURICCHIO201536,BOFFI20142151} for the well-posedness of the discrete
problem and the associated error estimates.
The ellipticity on the kernel is satisfied when $\beta_2 > \betaone > 0$, although the numerical findings reported in~\cite{BoffiGastaldi2024} suggest that
this assumption might be relaxed. Other stable alternatives are studied in~\cite{BoffiGastaldi2024}.

Since $\W_h \coloneqq V_h \times V_{2,h} \subset\W$ and $\Lambda_h\subset\Q$ are conforming subspaces, boundedness of the discrete operators follows automatically from the
continuous case. Hence, we only need to show the ellipticity on the discrete kernel and the discrete inf-sup condition.
Let $B_h$ be the discrete version of $B$, defined by $B_h(u_h,u_{2,h}) = \R u_h - u_{2,h}= u_h|_{\Omega_2} - u_{2,h}$. The discrete kernel is defined as follows
\begin{equation}
    \mathbb{K}_h = \{(v_h, v_{2,h}) \in \W_h : \langle \mu_h, v_h|_{\Omega_2} - v_{2,h} \rangle = 0 \text{ for all } \mu_h \in \Lambda_h\}.
\end{equation}

Let $P_h \colon L^2(\Omt) \to V_{2,h}$ be the standard $L^2$-orthogonal projection onto $V_{2,h}$. For
any $w \in L^2(\Omt)$, $P_h w \in V_{2,h}$ and is defined through the condition
\[
    (P_h w, z)_{\Omt} = (w, z)_{\Omt} \quad \forall z \in V_{2,h}.
\]It is well-known (cf.~\cite{Ern_Guermond-FEM-2004}) that $P_h$ satisfies, for $w \in V_2$,
\begin{equation}
    \norm{P_h w}_{0,\Omt} \le \norm{w}_{0,\Omt}. \label{eq:Ph_stability_L2}
\end{equation}
Moreover, thanks to the quasi-uniformity assumption on $\mathcal{T}_{2,h}$, the projection $P_h$ is
stable in the $H^1$-seminorm (cf.~\cite{CrouzeixThomee1987}), i.e. there exists a $C>0$ such that
\begin{equation}\label{eq:Ph_stability_H1}
    \seminorm{P_h w}_{1,\Omt} \le C \seminorm{w}_{1,\Omt} \quad \forall w \in V_2.
\end{equation}

We now show that $P_h$ is also stable in the weighted norm $\norm{\cdot}_\ast$ defined in \eqref{eq:imm-norm}.

\begin{lemma}\label{lem:Ph_weighted}
    Let $P_h:L^2(\Omt) \to V_{2,h}$ be the $L^2$-orthogonal projection onto $V_{2,h}$.
    Then, there exists a constant $C>0$ independent of discretization parameters such that
    \[
        \norm{P_h w}_\ast \le C \norm{w}_\ast \quad \forall w \in V_2.
    \]
\end{lemma}
\begin{proof}
    We apply the definition of norm $\norm{\cdot}_\ast$ to $P_h w$:
    \begin{align*}
        \norm{P_h w}_\ast^2 & \!=\! \gap \seminorm{P_h w}_{1,\Omt}^2 + \betaone \norm{P_h w}_{0,\Omt}^2
        \!\le\! \gap C^2 \seminorm{w}_{1,\Omt}^2 + \betaone \norm{w}_{0,\Omt}^2
        \!\le\! \max\{C^2, 1\} \norm{w}_\ast^2,
    \end{align*}
    where we used~\eqref{eq:Ph_stability_H1} and the $L^2$-stability of $P_h$. Taking square roots gives the
    thesis with constant $\sqrt{\max\{C ^2,1\}}$.
\end{proof}

In the next Theorem, we establish the discrete ellipticity on the kernel.
\begin{theorem}\label{thm:ellipticity-discrete}
    Let $V_h$, $V_{2,h}$, and $\Lambda_h$ be the finite element spaces defined above. Then, the ellipticity on
    the discrete kernel holds true, i.e. there exists a constant $C > 0$ independent of the discretization parameters such that
    \[
        a(w_h, w_h) \ge C\norm{w_h}_W^2 \quad \forall w_h \in \mathbb{K}_h.
    \]
\end{theorem}
\begin{proof}
    Let $w_h=(u_h,u_{2,h})\in\mathbb{K}_h$. By definition of $a(\cdot,\cdot)$, we have
    \[
        a(w_h, w_h) = \betaone \seminorm{u_h}_{1,\Om}^2 + \gap \seminorm{u_{2,h}}_{1,\Omt}^2.
    \]Since $w_h \in \mathbb{K}_h$, $u_{2,h}$ is the $L^2$-orthogonal projection of $u_h|_{\Omt} = \R u_h$ onto $V_{2,h}$. Hence, the stability
    of the projection $P_h$ and Lemma~\ref{lem:poincare} give
    \[
        \norm{u_{2,h}}_{0,\Omt} = \norm{P_h \R u_h}_{0, \Omt} \le \norm{\R u_h}_{0,\Omt} \le C_P \seminorm{u_h}_{1,\Om}.
    \]

    Using $\betaone \norm{u_{2,h}}_{0,\Omt}^2 \leq \betaone C_P^2 \seminorm{u_h}_{1,\Om}^2 \le C_P^2 a(w_h,w_h)$ in the definition of
    the $W$-norm gives $\norm{w_h}_W^2 = a(w_h, w_h) + \betaone \norm{u_{2,h}}_{0,\Omt}^2 \leq (1 + C_P^2) a(w_h, w_h)$, which rearranged yields
    \[
        a(w_h, w_h) \ge \frac{1}{1 + C_P^2} \norm{w_h}_W^2.
    \]
\end{proof}

To show the discrete inf-sup condition, we need to introduce the following two lemmas.
\begin{lemma}\label{lem:dual_ineq}
    Let $\norm{\cdot}_{V_2^{\ast}}$ be defined as in~\eqref{eqn:flat}.
    Then, for
    $\mu_h\in V_{2,h}$, there exists a constant $\widetilde C_{\mathrm{inv}}>0$ depending only on $\mathcal{T}_{2,h}$ and on its
    mesh regularity such that
    \begin{equation}\label{eq:inv-dual}
        \norm{\mu_h}_{0,\Omt}\le \widetilde C_{\mathrm{inv}}h_2^{-1}\norm{\mu_h}_{V_2^{\ast}} .
    \end{equation}
\end{lemma}

\begin{proof}
    We assume $0 \ne \mu_h \in V_{2,h}$. Testing the supremum at $q=\mu_h$ gives
    $\norm{\mu_h}_{V_2^{\ast}} \ge \norm{\mu_h}_{0,\Omt}^2/\norm{\mu_h}_{V_2}$. Using a standard inverse inequality, we have
    $$\norm{\mu_h}_{V_2}^2=\seminorm{\mu_h}_{1,\Omt}^2+\norm{\mu_h}_{0,\Omt}^2
        \le(C_{\mathrm{inv}}^2 h_2^{-2}+1)\norm{\mu_h}_{0,\Omt}^2.$$
    Since $h_2\le \operatorname{diam}(\Omt)$, %
    \[
        C_{\mathrm{inv}}^2 h_2^{-2}+1
        =C_{\mathrm{inv}}^2 h_2^{-2}\Bigl(1+C_{\mathrm{inv}}^{-2}h_2^2\Bigr)
        \le C_{\mathrm{inv}}^2 h_2^{-2}\Bigl(1+C_{\mathrm{inv}}^{-2}\operatorname{diam}(\Omt)^2\Bigr)
        = \widetilde C_{\mathrm{inv}}^2h_2^{-2},
    \]where $\widetilde C_{\mathrm{inv}}^2\coloneqq C_{\mathrm{inv}}^2\bigl(1+C_{\mathrm{inv}}^{-2}\operatorname{diam}(\Omt)^2\bigr)$ is a constant which depends only on
    $C_{\mathrm{inv}}$ and $\operatorname{diam}(\Omt)$, and not on $h_2$. Thus, $\norm{\mu_h}_{V_2}/\norm{\mu_h}_{0,\Omt}\le\widetilde C_{\mathrm{inv}}h_2^{-1}$, which gives the thesis.
\end{proof}

\begin{lemma}\label{lem:backgr_norm_stab}
    For $\mu_h \in \Lambda_h$, let $u_h = u_h(\mu_h) \in V_h$ be the unique solution of
    \begin{equation}\label{eq:uh}
        \beta_1 (\nabla u_h, \nabla v_h)_{\Omega} = \langle \mu_h, R v_h \rangle
        \qquad \forall v_h \in V_h,
    \end{equation}
    and let $C_\star := C_{\mathrm{app}} C_{\mathrm{reg}} \widetilde{C}_{\mathrm{inv}} C_E$,
    with $\widetilde{C}_{\mathrm{inv}}$ from Lemma~\ref{lem:dual_ineq}, $C_E$ from Lemma~\ref{lem:ext}, and
    $C_{\mathrm{app}}$, $C_{\mathrm{reg}}$ the interpolation and elliptic-regularity
    constants on $\Omega$. Then
    \[
        \beta_1 |u_h|_{1,\Omega}^2 \le \|\mu_h\|_{\Q,\mathrm{backgr}}^2
        \quad\text{and}\quad
        \beta_1 |u_h|_{1,\Omega}^2 \ge
        \big(1 - C_\star^2 (h/h_2)^2\big)\,\|\mu_h\|_{\Q,\mathrm{backgr}}^2 .
    \]
    In particular, if $h/h_2 \le \kappa/C_\star$ for some fixed $\kappa \in (0,1)$, then
    \[
        \beta_1 |u_h|_{1,\Omega}^2 \ge (1-\kappa^2)\,\|\mu_h\|_{\Q,\mathrm{backgr}}^2 .
    \]
\end{lemma}

\begin{proof}

    We start by proving the first inequality. Let $z=z(\mu_h)\in\Hoo$ be the continuous representative defined by~\eqref{eq:Bc-star},
    \begin{equation}\label{eq:z-cont}
        \betaone\inner{\nabla z}{\nabla v}_\Om=\dual{\mu_h}{\R v}\qquad\forall v\in\Hoo,
    \end{equation}
    so that $\norm{\mu_h}_{\Q,\mathrm{backgr}}^2=\betaone\seminorm{z}_{1,\Om}^2$ by
    \eqref{eq:Bc-extnorm}.

    Since $V_h\subset\Hoo$, Galerkin orthogonality follows after restricting~\eqref{eq:z-cont} to test functions in $V_h$. Therefore, we
    subtract~\eqref{eq:uh} from~\eqref{eq:z-cont} and obtain
    \begin{equation}\label{eq:galerkin}
        \betaone\inner{\nabla(z-u_h)}{\nabla v_h}_\Om=0\qquad\forall v_h\in V_h.
    \end{equation}

    Next, we decompose $z$ as $z=u_h+(z-u_h)$ and exploit bilinearity:
    \begin{equation*}
        \begin{aligned}
            \betaone \seminorm{z}_{1,\Om}^2
             & =\betaone \seminorm{u_h}_{1,\Om}^2
            +2 \betaone \inner{\nabla u_h}{\nabla(z-u_h)}_\Om
            +\betaone\seminorm{z-u_h}_{1,\Om}^2                                      \\
             & =\betaone\seminorm{u_h}_{1,\Om}^2+\betaone\seminorm{z-u_h}_{1,\Om}^2,
        \end{aligned}
    \end{equation*}
    the cross term vanishing by Galerkin orthogonality condition~\eqref{eq:galerkin}. Rearranging and
    using the fact that $\betaone \seminorm{z}_{1,\Om}^2=\norm{\mu_h}_{\Q,\mathrm{backgr}}^2$ gives
    \begin{equation}\label{eq:pyth}
        \betaone\seminorm{u_h}_{1,\Om}^2=\norm{\mu_h}_{\Q,\mathrm{backgr}}^2-\betaone\seminorm{z-u_h}_{1,\Om}^2,
    \end{equation}and the first assertion follows from $\betaone\seminorm{z-u_h}_{1,\Om}^2 \ge 0$.

    To obtain the second inequality, we bound the error term. Using Galerkin orthogonality
    in~\eqref{eq:galerkin} and writing $z-v_h=(z-u_h)+(u_h-v_h)$ yields
    \[
        \betaone \seminorm{z-v_h}_{1,\Om}^2=\betaone\seminorm{z-u_h}_{1,\Om}^2+\betaone\seminorm{u_h-v_h}_{1,\Om}^2\ \ge\ \betaone\seminorm{z-u_h}_{1,\Om}^2.
    \]
    Dividing by the common factor $\betaone$ and taking the minimum over $V_h$ gives
    \begin{equation}\label{eq:bestapprox}
        \seminorm{z-u_h}_{1,\Om}=\min_{v_h\in V_h}\seminorm{z-v_h}_{1,\Om}.
    \end{equation}

    Equation~\eqref{eq:z-cont} shows $z$ solves $-\betaone\Delta z=\mu_h\ind_{\Omt}$ in $\Om$. Since $\Omega$ is convex\footnote{As stated at the beginning of this section, $\Omega$ is a $d$-dimensional box.}, using elliptic
    regularity~\cite{grisvard2011elliptic} and the identity $\norm{\mu_h\ind_{\Omt}}_{0,\Om}=\norm{\mu_h}_{0,\Omt}$, we obtain
    \[
        \seminorm{z}_{2,\Om}\le C_{\mathrm{reg}}\betaone^{-1}\norm{\mu_h}_{0,\Omt},
    \]for a positive constant $C_{\mathrm{reg}}$ depending only on $\Om$. Let now $I_h$ denote the standard Lagrangian
    interpolant. Classical approximation estimates yield %
    \[
        \begin{aligned}
            \seminorm{z-u_h}_{1,\Om}=\min_{v_h\in V_h}\seminorm{z-v_h}_{1,\Om}
             & \le\seminorm{z-I_h z}_{1,\Om}\le C_{\mathrm{app}}h\seminorm{z}_{2,\Om}   \\
             & \le C_{\mathrm{app}}C_{\mathrm{reg}}\betaone^{-1}h\norm{\mu_h}_{0,\Omt},
        \end{aligned}
    \]where $C_{\mathrm{app}}$ is the constant associated with the error estimate for the Lagrangian
    interpolant $I_h$. We then apply
    the dual inequality in Lemma~\ref{lem:dual_ineq} and then Corollary~\ref{cor:Aext} $\bigl(\norm{\mu_h}_{V_2^{\ast}}\le C_E\betaone^{1/2}\norm{\mu_h}_{\Q,\mathrm{backgr}}\bigr)$, obtaining
    \[
        \seminorm{z-u_h}_{1,\Om}
        \le C_{\mathrm{app}}C_{\mathrm{reg}}\betaone^{-1}h\cdot \widetilde C_{\mathrm{inv}}h_2^{-1}
        \norm{\mu_h}_{V_2^{\ast}}
        \le C_\star\frac{h}{h_2}\betaone^{-1/2}\norm{\mu_h}_{\Q,\mathrm{backgr}}.
    \]
    Squaring and multiplying by $\betaone$, we get $\betaone\seminorm{z-u_h}_{1,\Om}^2\le C_\star^2(h/h_2)^2\norm{\mu_h}_{\Q,\mathrm{backgr}}^2$. Finally,
    inserting this into \eqref{eq:pyth} gives the second assertion.

    If $C_\star(h/h_2)\le\kappa$, for some fixed $\kappa\in(0,1)$, then $\betaone\seminorm{z-u_h}_{1,\Om}^2\le\kappa^2\norm{\mu_h}_{\Q,\mathrm{backgr}}^2$, and the last assertion follows.
\end{proof}

In the next Theorem, we establish the discrete inf-sup condition using the results developed in previous lemmas.
\begin{theorem}\label{thm:disc}
    Let $V_h, V_{2,h}, \Lambda_h$ be the equal-order finite element subspaces defined above
    for Problem~\ref{prob:discrete}, and let $\W_h \subset W$ and $\Lambda_h \subset \Lam$. If $h/h_2\le\kappa/C_\star$ for a fixed $\kappa\in(0,1)$ and
    the positive constant $C_\star$ defined as in Lemma~\ref{lem:backgr_norm_stab},
    then there exists a constant $\beta>0$ such that
    \[
        \inf_{0\ne \mu_h\in\Lambda_h} \sup_{0\ne w_h\in\W_h}\frac{\dual{B_h w_h}{\mu_h}}{\norm{w_h}_\W \norm{\mu_h}_\Q}
        \ \ge\ \beta.
    \]
\end{theorem}

\begin{proof}
    Let $0 \ne \mu_h\in\Lambda_h$. As in~\cite{prec1d2d}, we follow the continuous inf-sup in reverse order.
    We consider the discrete counterpart of the continuous case, i.e., we take
    \[
        w_h=(u_h,\,-P_h r)\in\W_h ,
    \]
    where $u_h = u_h(\mu_h) \in V_h$ is the solution defined in Lemma~\ref{lem:backgr_norm_stab} and $r=S^{-1}\mu_h
        \in V_2$ is the continuous representative as in Lemma~\ref{lem:dualnorm}. The numerator is
    then $$\dual{B_h w_h}{\mu_h}=\dual{\mu_h}{\R u_h}+\dual{\mu_h}{P_h r}.$$ From~\eqref{eq:uh} with $v_h=u_h$,
    $\dual{\mu_h}{\R u_h}=\betaone\seminorm{u_h}_{1,\Om}^2$. For the second term in the sum, we use the definition of $P_h$
    and the identity $(\mu_h,q)_{\Omt}=(r,q)_\ast$, obtaining
    \[
        \dual{\mu_h}{P_h r}=(\mu_h,P_h r)_{\Omt}=(\mu_h,r)_{\Omt}=(r,r)_\ast
        =\norm{\mu_h}_{\Q,\mathrm{imm}}^2 .
    \]
    Hence, the numerator can be expressed as $\dual{B_h w_h}{\mu_h}
        =\betaone\seminorm{u_h}_{1,\Om}^2+\norm{\mu_h}_{\Q,\mathrm{imm}}^2.$ Next, using the definition of $\norm{\cdot}_W$ and Lemma~\ref{lem:Ph_weighted}, we have
    \[
        \norm{w_h}_\W^2
        =\betaone\seminorm{u_h}_{1,\Om}^2+\norm{P_h r}_\ast^2
        \le \betaone\seminorm{u_h}_{1,\Om}^2+C^2\norm{\mu_h}_{\Q,\mathrm{imm}}^2 .
    \]
    Therefore,
    \[
        \frac{\dual{B_h w_h}{\mu_h}}{\norm{w_h}_\W}
        \ \ge\
        \frac{\betaone\seminorm{u_h}_{1,\Om}^2+\norm{\mu_h}_{\Q,\mathrm{imm}}^2}
        {\bigl(\betaone\seminorm{u_h}_{1,\Om}^2+C^2\norm{\mu_h}_{\Q,\mathrm{imm}}^2\bigr)^{1/2}}
        \ \ge\
        \frac{\bigl(\betaone\seminorm{u_h}_{1,\Om}^2+\norm{\mu_h}_{\Q,\mathrm{imm}}^2\bigr)^{1/2}}
        {\max\{1,C\}} .
    \]
    By Lemma~\ref{lem:backgr_norm_stab}, and since $1-\kappa^2\le1$,
    \[
        \betaone\seminorm{u_h}_{1,\Om}^2+\norm{\mu_h}_{\Q,\mathrm{imm}}^2
        \ \ge\ (1-\kappa^2)\norm{\mu_h}_{\Q,\mathrm{backgr}}^2+\norm{\mu_h}_{\Q,\mathrm{imm}}^2
        \ \ge\ (1-\kappa^2)\,\norm{\mu_h}_\Q^2 .
    \]
    The last two inequalities give $\frac{\dual{B_h w_h}{\mu_h}}{\norm{w_h}_\W}\ge\beta\norm{\mu_h}_\Q$, with
    $\beta=\sqrt{1-\kappa^2}/\max\{1,C\}$.
\end{proof}

\begin{remark}
    The result above shows that the inf-sup constant is independent of the mesh sizes, but \emph{not} of their ratio. Moreover, discretization spaces should be chosen so that their approximation properties are balanced both in terms of polynomial degrees, and in terms of the mesh ratio. Otherwise, the global error is dominated by the least accurate variable. Hence, the refinement of grids must be done simultaneously, and the proposed preconditioner is designed to be effective under these conditions.
\end{remark}

\subsection{Discrete linear system and preconditioner's realization}
We introduce the linear system associated with the discrete Problem~\ref{prob:discrete} and the corresponding preconditioner
matrix. Let $V_h=\text{span}\{\phi_i\}_{i=1}^{n}$, $V_{2,h}=\text{span}\{\psi_i\}_{i=1}^{m}$, and $\Lambda_h=\text{span}\{\chi_i\}_{i=1}^{p}$ be the finite element bases for the background, immersed, and multiplier spaces, respectively. Notice
that $\Lambda_h \equiv V_{2,h}$, and hence $p=m$ and $\{\chi_i\}_{i=1}^{p}=\{\psi_i\}_{i=1}^{m}$.
The matrices have sizes $\mat{A}_1 \in \mathbb{R}^{n \times n}$, $\mat{A}_2 \in \mathbb{R}^{m \times m}$,
$\mat{C} \in \mathbb{R}^{m \times n}$, and $\mat{M} \in \mathbb{R}^{m \times m}$.
We define their entries as follows:
\begin{equation}\label{eq:discrete_matrices}
    \begin{aligned}
        (\mat{A}_1)_{i,j} & = \inner{\nabla \phi_j}{\nabla \phi_i}_\Om ,
        \qquad            & (\mat{A}_2)_{i,j}                            & = \inner{\nabla \psi_j}{\nabla \psi_i}_{\Omt} , \\
        (\mat{C})_{i,j}   & = \inner{\chi_i}{\phi_j}_{\Omt} ,
        \qquad            & (\mat{M})_{i,j}                              & = \inner{\chi_j}{\chi_i}_{\Omt} .
    \end{aligned}
\end{equation}
With definitions in~\eqref{eq:discrete_matrices}, the matrix $\bbA$ of size $(n+m+p) \times (n+m+p)$ associated with
Problem~\ref{prob:discrete} is
\begin{equation}\label{eq:discrete_system}
    \bbA \coloneqq
    \begin{bmatrix}
        \betaone \mat{A}_1 & \mat{0}       & \mat{C}^\top  \\
        \mat{0}            & \gap\mat{A}_2 & -\mat{M}^\top \\
        \mat{C}            & -\mat{M}      & \mat{0}
    \end{bmatrix}.
\end{equation}
The coupling matrix $\mat{C}$ in~\eqref{eq:discrete_matrices} is the only non-trivial term to compute, as it involves products
of basis functions defined on the \emph{independent} meshes $\mathcal{T}_h$ and $\mathcal{T}_{2,h}$. Its assembly therefore
requires geometric search and collision detection strategies to identify overlapping background and
immersed elements. The analysis above assumes all integrals are computed exactly. In practice, the integrals defining $\mat{C}$ are approximated
by quadrature on immersed cells using tensor-product Gauss-Legendre rules with three points per coordinate direction. The background basis functions are evaluated at
these immersed quadrature points, without constructing mesh intersections. Related conditioning, stability, and
quadrature-error estimates are discussed in~\cite{BoffiConditioning,boffi2024quadratureerrorestimatesnonmatching}. More details on the assembly strategies of this matrix can be found in~\cite{boffi2023comparison,BOFFI2022115650}.

Finally, we shall represent the preconditioner $\Bcal$ by the matrix
\begin{equation}\label{eq:precond_discrete}
    \bbB \coloneqq \diag\left(
    (\betaone\mat{A}_1)^{-1},
    \bigl(\gap \mat{A}_2+\betaone \mat{M}\bigr)^{-1},
    \mat{M}^{-1}\mat{N}\mat{M}^{-1}
    \right),
\end{equation}
where $\mat{N} \coloneqq \theta \mat{A}_2 + \betaone \mat{M}$ is the matrix representation of the
operator $N$.
We emphasize that $\mat{N}$ alone is not a suitable representation for the multiplier block in the preconditioner $\bbB$. As discussed in~\cite[Section 4.1]{prec1d2d}, the inverse
of the mass matrix is used in~\eqref{eq:precond_discrete} as a natural adapter to obtain a matrix
operating between spaces suitable for preconditioning. We will make this precise in the following lemma.

\begin{lemma}\label{lem:discrete-multiplier}
    Let $\mathsf{N} \coloneqq \theta\,\mathsf{A}_2 + \beta_1\mathsf{M}$ be as in \eqref{eq:precond_discrete},
    $\mu_h\in \Lambda_h$ have coefficient vector $\mathbf{x}$, and $C$ be the
    constant in \eqref{eq:Ph_stability_H1}. Then
    \[
        C_E^{-2}\,
        \mathbf{x}^{\top}\mathsf{M}\mathsf{N}^{-1}\mathsf{M}\,\mathbf{x}
        \;\le\; \|\mu_h\|_{\Q}^2 \;\le\;
        2\,(1+C_P^2)\,\max\{C^2,1\}\,
        \mathbf{x}^{\top}\mathsf{M}\mathsf{N}^{-1}\mathsf{M}\,\mathbf{x},
    \]
    that is, $\mathsf{M}^{-1}\mathsf{N}\mathsf{M}^{-1}$ is spectrally equivalent,
    uniformly in discretization and model parameters, to the matrix representation of the Riesz map
    associated with the norm $\|\cdot\|_{\Lam}$ on $\Lambda_h$.
\end{lemma}

\begin{proof}
    The case $\mu_h =0$ is trivial, so we assume $\mu_h \neq 0$. Recall that $N=\theta(-\Delta^N_{\Omega_2})+\beta_1 I$, and $\mat{N}$ denotes its matrix representation. Since $\theta$ and $\beta_1$ are positive, $\mat{N}$ is an SPD matrix. We introduce on $V_2$ the inner product $(q,p)_N \coloneqq \theta(\nabla q,\nabla p)_{\Omega_2} + \beta_1(q,p)_{\Omega_2}$,
    with induced norm $\|\cdot\|_N$.  Let $r_h\in V_{2,h}$ be the solution of
    $$(r_h,q_h)_N=(\mu_h,q_h)_{\Omega_2} \quad \forall q_h\in V_{2,h},$$ which in matrix form reads
    $\mathsf{N}\mathbf{r}=\mathsf{M}\mathbf{x}$. Then, the following chain of equalities holds:
    \begin{equation}\label{eq:xMNx}
        \begin{aligned}
            \mathbf{x}^{\top}\mathsf{M}\mathsf{N}^{-1}\mathsf{M}\mathbf{x}
             & = \mathbf{x}^{\top}\mathsf{M}\mathbf{r}
            = (\mu_h,r_h)_{\Omega_2}
            = \|r_h\|_N^2                              \\
             & = \sup_{0\neq q_h\in V_{2,h}}
            \frac{(\mu_h,q_h)_{\Omega_2}^2}{\|q_h\|_N^2}
            \eqqcolon \|\mu_h\|_{N,\ast,h}^2.
        \end{aligned}
    \end{equation}
    The last two identities follow as in Lemma~\ref{lem:dualnorm}. Again following the argument of
    Lemma~\ref{lem:dualnorm}, and noting that the pairings coincide with the $L^2(\Omega_2)$ inner product since $\mu_h \in \Lambda_h \subset L^2(\Omega_2)$, we set
    \[
        \|\mu_h\|_{N,\ast}^2 \coloneqq \sup_{0\neq q\in V_2}
        \frac{(\mu_h,q)_{\Omega_2}^2}{\|q\|_N^2}
        =(\mu_h,\,N^{-1}\mu_h)_{\Omega_2}.
    \]Since
    $V_{2,h}\subset V_2$, it follows that $\|\mu_h\|_{N,\ast,h}\le\|\mu_h\|_{N,\ast}$.

    Conversely, let $0\neq q\in V_2$, and $P_h$ be the $L^2(\Omega_2)$-orthogonal projection onto $V_{2,h}$.
    Using the definition of $P_h$, testing the supremum in~\eqref{eq:xMNx} at $q_h=P_h q$, and Lemma~\ref{lem:Ph_weighted}, we have\footnote{Lemma~\ref{lem:Ph_weighted} holds verbatim also for $\|\cdot\|_N$ with the same constant
        $\max\{C,1\}$, since its proof uses only~\eqref{eq:Ph_stability_L2} and~\eqref{eq:Ph_stability_H1}, which do not depend on the weights.}:
    \[
        (\mu_h,q)_{\Omega_2}
        =(\mu_h,P_hq)_{\Omega_2}
        \le \|\mu_h\|_{N,\ast,h}\,\|P_hq\|_N
        \le \max\{C,1\}\,\|\mu_h\|_{N,\ast,h}\,\|q\|_N.
    \]
    Dividing by $\|q\|_N$, and taking the supremum over $0\neq q\in V_2$
    yields $\|\mu_h\|_{N,\ast}\le\max\{C,1\}\,\|\mu_h\|_{N,\ast,h}$. Therefore
    \begin{equation}\label{eq:Nh_bound}
        \|\mu_h\|_{N,\ast,h}^2 \;\le\; (\mu_h,\,N^{-1}\mu_h)_{\Omega_2} \;\le\;
        \max\{C^2,1\}\,\|\mu_h\|_{N,\ast,h}^2 .
    \end{equation}

    Finally, the equivalence in Theorem~\ref{thm:multiplier_block} gives, for all $\mu_h\in \Lambda_h\subset\Lambda$,
    \begin{equation}\label{eq:multiplier_block_bound}
        C_E^{-2}\,(\mu_h,\,N^{-1}\mu_h)_{\Omega_2}
        \le \|\mu_h\|_{\Q}^2
        \le 2(1+C_P^2)\,(\mu_h,\,N^{-1}\mu_h)_{\Omega_2}.
    \end{equation}
    Combining~\eqref{eq:multiplier_block_bound} with \eqref{eq:Nh_bound} and \eqref{eq:xMNx} yields the thesis.
\end{proof}

Theorems~\ref{thm:ellipticity-discrete} and \ref{thm:disc} with Lemma~\ref{lem:discrete-multiplier}, together with the operator preconditioning
framework of~\cite{MW}, yield a spectral condition number $\kappa(\bbB \bbA)$ bounded independently of the discretization and of the
jump, provided the mesh ratio condition is satisfied.

Having in mind large scale three-dimensional problems, exact solves for the diagonal blocks of the preconditioner $\bbB$ in~\eqref{eq:precond_discrete} are not feasible. Hence, we will
also employ an approximate variant of the preconditioner $\bbB$, denoted by $\bbBt$, where the first two diagonal blocks are inverted approximately using a fixed number of V-cycles of algebraic (or geometric) multigrid.

\section{Numerical Experiments}\label{sec:experiments}

We perform numerical experiments to validate the theoretical results
and assess the practical performance of the proposed preconditioner. Computations have been carried out on a laptop equipped with an AMD Ryzen 9 PRO 7940HS processor and 64 GB of RAM. Experiments
have been performed using the \textsc{C++} finite element library \textsc{deal.II}~\cite{dealIIdesign}, which provides efficient handling of non-matching meshes in a parallel setting, as well as
algebraic multigrid (AMG) preconditioners exposed by well-established linear algebra libraries. The source code implementing the methodology and the experiments presented is
available at the dedicated GitHub repository provided by the author.

The outer solver is MinRes~\cite{PaigeSaunders1975}. We first validate the mesh independence and robustness of the preconditioner $\bbB$ (and of its more practical variant $\bbBt$), and then we compare them with the
modified augmented Lagrangian preconditioner proposed in recent joint work~\cite{BENZI2026522}%
\footnote{In~\cite{BENZI2026522}, AL preconditioners are compared with block-triangular preconditioning from~\cite{alshehri2025multigridpreconditioningfddlmmethod}, showing bounded iteration counts
    for large coefficient jumps.}. The evaluation metrics are the robustness against coefficient jumps, the number
of iterations required for convergence, and the associated wall-clock time.

As usual in fictitious domain methods, the background domain $\Omega$ is discretized with a structured mesh $\mathcal{T}_h$ consisting of
quadrilaterals or hexahedra. In all experiments, we report the ratio $h/h_2$ in each caption, and study refinement at fixed ratio\footnote{In the experiments below, $h/h_2$ ranges from about $0.2$ to $0.7$, delivering robust performance.}. All iterations start from a zero initial guess and, unless otherwise specified, are stopped
when the relative residual norm drops below $10^{-6}$, or when the absolute residual is below $10^{-9}$.

\subsection{Validation of the preconditioner}\label{subsec:validation_2d}
We start with a two-dimensional experiment, where the background and immersed domains are set
to $\Omega=(-1,1)^2$ and $\Omega_2 = (-0.24,0.44)^2$, respectively. To ensure fully non-matching grids, we rotate the mesh
for $\Omega_2$ by $30$ degrees. We investigate the robustness of the proposed preconditioner with respect to the contrast between $\betaone$ and $\rtwo$ as well as
against mesh refinement. We set $\betaone = 1$ and vary $\rtwo$ to increase the contrast between subdomains. The diagonal blocks of the preconditioner $\bbBt$ are inverted using three V-cycles of the AMG preconditioner from the \textsc{TrilinosML} package~\cite{Trilinos} (version 14.4.0). We use a Chebyshev smoother with three sweeps,
\textsc{Amesos-KLU} as the coarse solver, an aggregation threshold of $10^{-2}$, and a maximum coarse-level size of 2000. Mass matrices appearing
as space adapters in the multiplier block are inverted exactly with \textsc{MUMPS}~\cite{MUMPS}, since their size is the number of degrees of freedom (dofs) for the multiplier space, which is relatively
small in these examples.
Table~\ref{tab:2d_square} reports the number of MinRes iterations across mesh refinement levels for increasing values of the jump $\alpha$. Iterations are bounded
and stable with respect to both mesh refinement and parameters.

To test the sensitivity of the preconditioner to the geometry of the immersed domain, we repeat
the same experiment with the circular immersed domain centered at the origin with radius $0.3$, i.e., $\Omega_2 = \mathcal{B}_{0.3}(0,0)$, and report the results in
Table~\ref{tab:2d_ball}, confirming the findings observed in the previous experiment.

In parentheses in Table~\ref{tab:2d_square} and Table~\ref{tab:2d_ball}, we show the effect of removing the weights $\alpha$, $\betaone$, and $\theta$, i.e., setting $\alpha = \betaone = \theta = 1$, from the preconditioner $\bbB$ in~\eqref{eq:precond_discrete} on the convergence behavior of the iterative solver. To isolate the
effect of dropping the weights from other potential sources of iteration growth, we invert the diagonal blocks with the direct solver \textsc{MUMPS}. For each fixed row, the iterations are evidently increasing with the jump, demonstrating
the loss of robustness with respect to the coefficient contrast of the preconditioner without the predicted weights.

To further assess robustness, Table~\ref{tab:kappa-2d-square} reports estimates of the spectral condition number $\kappa (\bbB \bbA)$ under mesh refinements and varying
coefficient contrasts for the immersed square test case. Matrices have been exported for eigenvalue analysis in Python, and $\kappa(\bbB \bbA)$ has been estimated as the ratio between the largest
and smallest eigenvalues in modulus of the generalized eigenproblem $\bbA x = \lambda \bbB^{-1} x$ using the \texttt{scipy.sparse.linalg.eigsh} function.
The estimates remain bounded across both coefficient contrasts and mesh refinements, with $\kappa_{\min}=5.120$, $\kappa_{\max}=7.524$.

\begin{table}[ht]
    \centering
    \begin{minipage}[t]{0.49\textwidth}
        \vspace{0pt}
        \centering
        \sisetup{group-separator={,}, group-minimum-digits=4}
        \resizebox{\linewidth}{!}{%
            \begin{tabular}{
                    S[table-format=7]
                    S[table-format=7]
                    S[table-format=6]
                    S[table-format=2]
                    S[table-format=2]
                    S[table-format=2]
                    S[table-format=2]
                }
                \toprule
                \multicolumn{7}{c}{\textbf{FD-DLM (2D, immersed square)}}                                                                                                               \\
                \midrule
                {$\mathrm{dof}_u$} & {$\mathrm{dof}_{u_2}$} & {$\mathrm{dof}_\lambda$} & \multicolumn{4}{c}{MinRes iterations}                                                          \\
                \cmidrule(lr){4-7}
                                   &                        &                          & {$\beta_2=10$}                        & {$\beta_2=10^3$} & {$\beta_2=10^5$} & {$\beta_2=10^7$} \\
                \midrule
                4225               & 81                     & 81                       & {23 (31)}                             & {26 (51)}        & {26 (66)}        & {26 (86)}        \\
                16641              & 289                    & 289                      & {24 (33)}                             & {26 (54)}        & {26 (70)}        & {27 (94)}        \\
                66049              & 1089                   & 1089                     & {25 (34)}                             & {28 (56)}        & {28 (74)}        & {28 (95)}        \\
                263169             & 4225                   & 4225                     & {25 (34)}                             & {29 (56)}        & {30 (76)}        & {30 (95)}        \\
                1050625            & 16641                  & 16641                    & {27 (36)}                             & {29 (58)}        & {30 (77)}        & {30 (96)}        \\
                4198401            & 66049                  & 66049                    & {26 (36)}                             & {28 (58)}        & {29 (77)}        & {29 (96)}        \\
                \bottomrule
            \end{tabular}}
        \caption{Immersed square: MinRes iteration counts for $\bbBt$ across mesh refinements (ratio $h/h_2=0.37$) and coefficient jumps. Between parentheses the iteration counts for the preconditioner $\bbB$ \emph{without} the
            weights $\alpha$, $\betaone$, and $\theta$.}
        \label{tab:2d_square}
    \end{minipage}\hfill
    \begin{minipage}[t]{0.49\textwidth}
        \vspace{0pt}
        \centering
        \sisetup{group-separator={,}, group-minimum-digits=4}
        \resizebox{\linewidth}{!}{%
            \begin{tabular}{
                    S[table-format=7]
                    S[table-format=7]
                    S[table-format=6]
                    S[table-format=2]
                    S[table-format=2]
                    S[table-format=2]
                    S[table-format=2]
                }
                \toprule
                \multicolumn{7}{c}{\textbf{FD-DLM (2D, immersed ball)}}                                                                                                                 \\
                \midrule
                {$\mathrm{dof}_u$} & {$\mathrm{dof}_{u_2}$} & {$\mathrm{dof}_\lambda$} & \multicolumn{4}{c}{MinRes iterations}                                                          \\
                \cmidrule(lr){4-7}
                                   &                        &                          & {$\beta_2=10$}                        & {$\beta_2=10^3$} & {$\beta_2=10^5$} & {$\beta_2=10^7$} \\
                \midrule
                4225               & 17                     & 17                       & {15 (21)}                             & {15 (28)}        & {15 (36)}        & {15 (44)}        \\
                16641              & 57                     & 57                       & {17 (22)}                             & {18 (35)}        & {18 (47)}        & {18 (58)}        \\
                66049              & 209                    & 209                      & {19 (25)}                             & {20 (40)}        & {20 (53)}        & {20 (66)}        \\
                263169             & 801                    & 801                      & {21 (28)}                             & {21 (43)}        & {21 (56)}        & {21 (70)}        \\
                1050625            & 3137                   & 3137                     & {19 (31)}                             & {21 (55)}        & {21 (76)}        & {21 (97)}        \\
                4198401            & 12417                  & 12417                    & {19 (29)}                             & {19 (50)}        & {18 (69)}        & {18 (99)}        \\
                \bottomrule
            \end{tabular}}
        \caption{Immersed ball: MinRes iteration counts for $\bbBt$ across mesh refinements (ratio $h/h_2=0.19$) and coefficient jumps. Between parentheses the iteration counts for the preconditioner $\bbB$ \emph{without} the
            weights $\alpha$, $\betaone$, and $\theta$.}
        \label{tab:2d_ball}
    \end{minipage}
\end{table}

\begin{table}[ht]
    \centering
    \small
    \begingroup
    \footnotesize
    \sisetup{group-separator={,}, group-minimum-digits=4}
    \setlength{\tabcolsep}{3pt}
    \renewcommand{\arraystretch}{1.0}
    \begin{tabular}{
            @{}
            S[table-alignment-mode=none,table-number-alignment=right,group-digits=false,output-exponent-marker={}]
            S[table-alignment-mode=none,table-number-alignment=right]
            S[table-format=1.3]
            S[table-format=1.3]
            S[table-format=1.3]
            S[table-format=1.3]
            S[table-format=1.3]
            S[table-format=1.3]
            @{}
        }
        \toprule
        \multicolumn{8}{c}{\textbf{FD-DLM (2D, immersed square). Condition numbers $\kappa(\bbB\bbA)$}}                                                  \\
        \midrule
        {$h$}       & {$N_{\text{dof}}$} & {$\beta_2=10$} & {$\beta_2=10^3$} & {$\beta_2=10^4$} & {$\beta_2=10^5$} & {$\beta_2=10^6$} & {$\beta_2=10^7$} \\
        \midrule
        6.25e-2     & 1139               & 5.120          & 5.853            & 5.860            & 5.861            & 5.861            & 5.861            \\
        3.125e-2    & 4387               & 5.487          & 6.406            & 6.416            & 6.417            & 6.417            & 6.417            \\
        1.5625e-2   & 17219              & 5.754          & 6.824            & 6.836            & 6.837            & 6.837            & 6.837            \\
        7.8125e-3   & 68227              & 5.939          & 7.120            & 7.133            & 7.134            & 7.134            & 7.134            \\
        3.90625e-3  & 271619             & 6.072          & 7.338            & 7.352            & 7.353            & 7.353            & 7.353            \\
        1.953125e-3 & 1083907            & 6.175          & 7.508            & 7.523            & 7.524            & 7.524            & 7.524            \\
        \bottomrule
    \end{tabular}
    \endgroup
    \caption{Condition numbers $\kappa(\bbB\bbA)$ for the immersed-square setup over mesh refinements and coefficient contrasts. $h$ denotes the background mesh size (ratio $h/h_2=0.37$), and $N_{\text{dof}}$ the number of total degrees
        of freedom. The parameter $\betaone$ is set to $1$.}
    \label{tab:kappa-2d-square}
\end{table}

\subsection{Comparison with AL preconditioner from~\texorpdfstring{\cite{BENZI2026522}}{}}
\begin{figure}[!ht]
    \centering
    \begin{minipage}[t]{0.49\textwidth}
        \centering
        \subcaptionbox{$\Omega_2=(-0.24,0.44)^2$ (rotated by $30^\circ$).\label{fig:time-comparison-al-vs-b}}[\linewidth]{%
            \includegraphics[width=\linewidth]{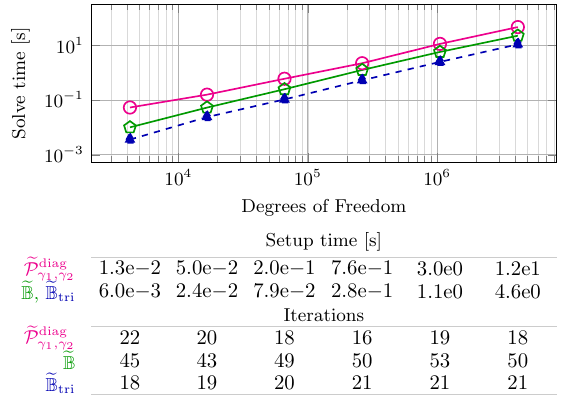}}
    \end{minipage}\hfill
    \begin{minipage}[t]{0.47\textwidth}
        \centering
        \subcaptionbox{$\Omega_2=(-0.24,0.44)^3$ (rotated by $30^\circ$).\label{fig:time-comparison-al-vs-b-3D}}[\linewidth]{%
            \includegraphics[width=\linewidth]{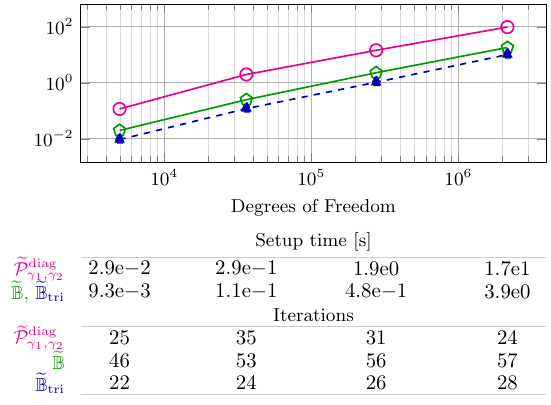}}
    \end{minipage}
    \caption{Wall-clock time comparisons (left: 2D, right: 3D) between the modified AL preconditioner (with parameters $\gamma_1=10$ and $\gamma_2=10^{-2}$) $\mathcal{\widetilde P}_{\gamma_1,\gamma_2}^{\mathrm{diag}}$~(\methodkey{o}) from~\cite{BENZI2026522}, the
    preconditioner $\bbBt$~(\methodkey{pentagon}), and its block-triangular variant $\bbBt_{\mathrm{tri}}$~(\methodkey[dashed]{triangle*}). Setup times [s] and iteration counts are listed by method below each plot, aligned with the corresponding problem sizes. The mesh ratio is as in Table~\ref{tab:2d_square}.}
    \label{fig:time-comparison-al-vs-b-combined}
\end{figure}

We compare the performance of $\bbB$ with the modified AL preconditioner $\mathcal{\widetilde P}_{\gamma_1,\gamma_2}^{\mathrm{diag}}$ recently developed in~\cite{BENZI2026522}.
The AL preconditioner $\mathcal{\widetilde P}_{\gamma_1,\gamma_2}^{\mathrm{diag}}$ employs FGMRES(30) as outer solver, and diagonal blocks
are approximately inverted with the conjugate gradient method, preconditioned by one V-cycle of AMG, with the same settings as before, and using a loose
absolute tolerance set to $10^{-2}$ for the inner CG iterations. Following~\cite{BENZI2026522}, parameters $\gamma_1$ and $\gamma_2$ are set to $10$ and $10^{-2}$, respectively. To ensure a fair comparison of the two preconditioners, we terminate the iterations when the norm of the
residual of the system falls below $10^{-9}$. Both preconditioners are run on the same machine and with the same software stack.
The setup cost of the preconditioners includes the construction of the AMG hierarchies for the diagonal blocks, the inversion of mass matrices, and, in the AL case, the construction of the augmented terms. The first two diagonal blocks in $\bbBt$ are approximately
inverted with one V-cycle, using the same AMG setup described above, while mass matrices in the multiplier block of $\bbBt$ are inverted exactly with \textsc{MUMPS}. We consider the previous two-dimensional
experiment with the immersed square, with $\beta_2=10^3$ and the same mesh ratio $h/h_2$ and number of degrees of freedom, and report the solver times for both strategies. Setup costs and iteration
counts are in aligned panels below the timing plots. Results in Figure~\ref{fig:time-comparison-al-vs-b} show that the
preconditioner $\bbBt$ is faster than the modified AL preconditioner, with a speedup of approximately $2\times$ for the largest problem size considered. Identical trends are observed for other contrasts, and hence we omit
reporting them here.

Likewise, we repeat the same experiment in three dimensions, with the background domain $\Omega=(-1,1)^3$ and the immersed domain $\Omega_2 = (-0.24,0.44)^3$. The wall-clock
times are reported in Figure~\ref{fig:time-comparison-al-vs-b-3D}. The gap between the two approaches is evident, with a clear gain in favor of $\bbBt$: for the finest discretization, it achieves
convergence in approximately $18$ seconds, while the modified AL preconditioner $\mathcal{\widetilde P}_{\gamma_1,\gamma_2}^{\mathrm{diag}}$ requires
around $97$ seconds, resulting in a speedup of more than $5\times$.
These results show that $\bbBt$ outperforms the modified AL preconditioner in terms of wall-clock time, while maintaining robustness with respect
to the contrast between coefficients and mesh refinement.

In both cases, $\bbB$ and the practical variant $\bbBt$ exhibit higher iteration counts compared to the AL competitor $\mathcal{\widetilde P}_{\gamma_1,\gamma_2}^{\mathrm{diag}}$. Nevertheless, their application
is much cheaper due to the simpler structure of the diagonal blocks, which can be efficiently inverted with standard and optimal multigrid solvers. In contrast, the AL preconditioner
requires complex and tailored solvers for the diagonal blocks and involves additional parameters that need to be tuned for optimal performance. While for linearizations of Navier--Stokes problems
such solvers have been developed over the years, see e.g.~\cite{ALprec,Farrell2019}, we stress that multilevel solvers for the augmented systems in fictitious domain discretizations are non-trivial and yet
to be developed. In addition, the AL preconditioner requires additional setup work to construct the augmented terms. These terms are absent
from $\bbBt$, explaining its substantially lower setup cost.
Both panels include also the following natural block-triangular variant $\bbBt_{\mathrm{tri}}$, used with GMRES\footnote{The block-triangular structure of $\bbBt_{\mathrm{tri}}$ makes the preconditioned system non-symmetric. We therefore use GMRES, with right preconditioning.}, defined by
the inverse of the following block-triangular matrix:
\begin{equation*}
    \begin{bmatrix}
        \betaone \mat{A}_1 &                                 & \mat{C}^\top                \\[1ex]
                           & \gap \mat{A}_2+\betaone \mat{M} & -\mat{M}                    \\[1ex]
                           &                                 & -\mat{M}\mat{N}^{-1}\mat{M}
    \end{bmatrix},
\end{equation*}where the first two diagonal blocks are inverted again using one V-cycle of AMG. It halves iteration counts (and correspondingly reduces wall-clock times) compared to the block-diagonal version, while
retaining the same computational gains.

We point out that $\bbBt$ is also much easier to implement compared to AL approaches: it does not require parameter tuning\footnote{Some sensitivity with respect to the jump's size $\beta_2-\beta_1$ has indeed been observed in~\cite{BENZI2026522}, and mitigated by adjusting values of $\gamma_1$ and $\gamma_2$.} for $\gamma_1$ and $\gamma_2$, it avoids the
need for tailored geometric multigrid solvers for the diagonal blocks, and does not modify the sparsity pattern of the original system matrix through
the augmented terms. In addition, the parallel implementation of $\bbBt$ is fairly straightforward.
In light of these reasons, we claim that $\bbBt$ is a more practical preconditioner for this problem, and we will
not pursue other comparisons with AL preconditioners in the subsequent experiments.

\subsection{Extension to immersed linear elasticity}
The analysis in Sections~\ref{sec:operator_prec}--\ref{sec:discretization} concerns the scalar, equal-order, discretization. We now test its
formal analogue for the immersed linear elasticity case. Let $\mu,\mu_2$ be the shear
moduli in $\Omega$ and $\Omega_2$, respectively, and $\lambda, \lambda_2$ the first Lamé parameters. Consistently with
the assumptions presented in Section~\ref{sec:model_problem}, we assume the coefficients to satisfy $\mu_2 > \mu>0$
and $\lambda_2 > \lambda>0$. We denote by $\delta_\mu = \mu_2 - \mu$ and $\delta_\lambda = \lambda_2 - \lambda$ the corresponding jumps. Vector-valued functions are indicated in boldface\footnote{Thus, no ambiguity arises when $\lambda$ denotes either a multiplier or
    a Lamé parameter.}. The
interested reader is referred to~\cite{ALZETTA2020106334} for the study of fictitious domain formulations applied to fiber-reinforced and multiscale materials. Augmented Lagrangian preconditioners for this problem have been recently
applied in~\cite{BENZI2026522,benzi2026unifiedvariationalframeworkaugmented}, and extended to beam--solid coupling in~\cite{FIRMBACH2027119425}. We fix the background parameters $(\mu,\lambda)=(1,2)$ and consider contrasts
$\mu_2/\mu=\lambda_2/\lambda\in\{10,10^2,10^3,10^4\}$.

The corresponding vector-valued formulation reads as follows.
Let $d \in\{2,3\}$, and let $\Gamma_{\mathrm{D}} \subset \partial \Omega$ be the Dirichlet boundary. Given $\boldsymbol{f} \in \bigl[L^2(\Omega)\bigr]^d$, find $(\boldsymbol{u},\boldsymbol{u}_2,\boldsymbol{\lambda})\in V\times V_2\times\Lambda$ such that
\begin{align*}\label{prob:elasticity}
    2 \mu(\varepsilon(\boldsymbol{u}),\varepsilon(\boldsymbol{v}))_\Omega + \lambda(\operatorname{div} \boldsymbol{u}, \operatorname{div} \boldsymbol{v})_\Omega + c(\boldsymbol{\lambda},\boldsymbol{v}|_{\Omega_2})                        & =(\boldsymbol{f},\boldsymbol{v})_\Omega &  & \forall \boldsymbol{v} \in V, \nonumber       \\
    2 \delta_\mu (\varepsilon(\boldsymbol{u}_2),\varepsilon(\boldsymbol{v}_2))_{\Omega_2} + \delta_\lambda(\operatorname{div} \boldsymbol{u}_2, \operatorname{div} \boldsymbol{v}_2)_{\Omega_2} -c(\boldsymbol{\lambda}, \boldsymbol{v}_{2}) & = 0                                     &  & \forall \boldsymbol{v}_{2}\in V_{2},          \\
    c(\boldsymbol{\mu},\boldsymbol{u}|_{\Omega_2} - \boldsymbol{u_{2}})                                                                                                                                                                      & = 0                                     &  & \forall \boldsymbol{\mu}\in\Lambda, \nonumber
\end{align*}where $\varepsilon(\boldsymbol{u}) = \frac{1}{2}(\nabla \boldsymbol{u} + \nabla \boldsymbol{u}^T)$ is the linearized strain tensor, $V = [H_{\Gamma_{\mathrm{D}}}^1(\Omega)]^d, V_2=[H^1(\Omega_2)]^d, \Lambda=\bigl([H^1(\Omega_2)]^d\bigr)^\ast$ are the vector-valued analogues
of the spaces defined for the scalar problem, and $c(\boldsymbol{\lambda},\boldsymbol{v}) = \langle \boldsymbol{\lambda},\boldsymbol{v} \rangle$.
Since $\mu_2/\mu=\lambda_2/\lambda$, the background and immersed
elasticity forms are multiples of a common normalized elasticity form,
with factors $2\mu+\lambda$ and $2\delta_\mu+\delta_\lambda$,
respectively. This gives the correspondence
$\betaone \leftrightarrow 2\mu+\lambda$ and
$\gap \leftrightarrow 2\delta_\mu+\delta_\lambda$
in the scalar construction.
Let $\mat{A}_{1,\mathrm{el}}(\mu,\lambda)$ be the stiffness matrix
of the elasticity form on $\Omega$ with Lamé parameters $\mu$
and $\lambda$, $\mat{A}_{2,\mathrm{el}}(\cdot,\cdot)$ its counterpart
on $\Omega_2$, and $\mat{M}$ the mass matrix on $\Lambda_h$.
The direct analogue of the preconditioner~\eqref{eq:precond_discrete} is then the
candidate preconditioner
\begin{equation}\label{eq:precond_el}
    \bbB_{\mathrm{el}} \coloneqq \diag\left(
    \mat{A}_{1,\mathrm{el}}(\mu,\lambda)^{-1},
    \bigl(
    \mat{A}_{2,\mathrm{el}}(\delta_\mu,\delta_\lambda)
    +(2\mu+\lambda)\mat{M}
    \bigr)^{-1},
    \mat{M}^{-1}\mat{N}_{\mathrm{el}}\mat{M}^{-1}
    \right),
\end{equation}
with
$\mat{N}_{\mathrm{el}}\coloneqq
    \mat{A}_{2,\mathrm{el}}(\theta_\mu,\theta_\lambda)
    +(2\mu+\lambda)\mat{M}$, where
$\theta_\mu\coloneqq\mu\delta_\mu/\mu_2$ and
$\theta_\lambda\coloneqq\lambda\delta_\lambda/\lambda_2$.

We start with a two-dimensional test with $\Omega=(-1,1)^2$, and $\Omega_2=(-0.24,0.44)^2$. As done in Section~\ref{subsec:validation_2d}, the mesh
for $\Omega_2$ is rotated by $30$ degrees to ensure that background and immersed grids are non-matching. We fix the boundary conditions to be homogeneous Dirichlet on the left side of $\Omega$, and homogeneous Neumann on the remaining sides. We apply
a body force in the negative y direction, defined as $\boldsymbol{f}(\boldsymbol{x}) = -1 \boldsymbol{e}_2 \ind_{\{x_1 > 0.6\}}(\boldsymbol{x})$. For each set of parameters, the materials exhibit a huge contrast in their Young's modulus $E=2\mu(1+\nu)$, which measures the material's stiffness under tensile or compressive force.
Notice that these numerical tests do not explore the nearly incompressible limit $\nu\to 0.5$, which is an interesting
direction for future investigation.

To assess the validity of the preconditioner in this two-dimensional experiment, we invert the diagonal blocks
in $\bbB_{\mathrm{el}}$ through the sparse direct solver \textsc{MUMPS}. Table~\ref{tab:2d_square_linear_elasticity} shows a sweep over different contrast values
between the Lamé parameters in the two subdomains, reporting the number of MinRes iterations under refinement.
Iteration counts remain bounded across mesh refinement, and stabilize to constant values. Therefore, in these tests the preconditioner $\bbB_{\mathrm{el}}$ is robust with respect to both
mesh refinement and contrast between the Lamé parameters, confirming the trends observed in the
scalar case. Figure~\ref{fig:displacement_square} shows the displacement field $\boldsymbol{u}_h$ for the case with $\mu_2/\mu=\lambda_2/\lambda=10^2$.

We repeat the same experiment in three dimensions, and denote by $\bbBt_{\mathrm{el}}$ the approximate variant of the preconditioner $\bbB_{\mathrm{el}}$. The same boundary conditions as in the two-dimensional case are applied, the immersed cube is rotated by $30$ degrees around the y-axis, and we employ
a body force in the negative z direction defined as $\boldsymbol{f}(\boldsymbol{x}) = -1 \boldsymbol{e}_3 \ind_{\{x_1 > 0.6\}}(\boldsymbol{x})$. We use a fixed number of five AMG V-cycles for the inversion of the diagonal blocks in $\bbBt_{\mathrm{el}}$. The six rigid-body modes are supplied as near null-space
vectors in the setup of the AMG hierarchies for the first two diagonal blocks. Table~\ref{tab:3d_cube_linear_elasticity} and Figure~\ref{fig:displacement_cube} show the results for this configuration, confirming
the robustness of the preconditioner with respect to both mesh refinement and contrast between the subdomains also in three dimensions and with approximate inversions of the diagonal blocks.

\begingroup
\setlength{\intextsep}{6pt}
\captionsetup{belowskip=0pt}
\begin{table}[ht]
    \centering
    \begin{minipage}[t]{0.49\textwidth}
        \vspace{0pt}
        \centering
        \sisetup{group-separator={,}, group-minimum-digits=4}
        \resizebox{\linewidth}{!}{%
            \begin{tabular}{
                    S[table-format=7]
                    S[table-format=7]
                    S[table-format=6]
                    S[table-format=2]
                    S[table-format=2]
                    S[table-format=2]
                    S[table-format=2]
                }
                \toprule
                \multicolumn{7}{c}{\textbf{Elasticity test case (2D, immersed square)}}                                                                                                                                \\
                \midrule
                {$\mathrm{dof}_{\boldsymbol{u}}$} & {$\mathrm{dof}_{\boldsymbol{u_2}}$} & {$\mathrm{dof}_{\boldsymbol{\lambda}}$} & \multicolumn{4}{c}{MinRes iterations}                                              \\
                                                  &                                     &                                         & \multicolumn{4}{c}{$\mu_2/\mu=\lambda_2/\lambda$}                                  \\
                \cmidrule(lr){4-7}
                                                  &                                     &                                         & {$10$}                                            & {$10^2$} & {$10^3$} & {$10^4$} \\
                \midrule
                2178                              & 162                                 & 162                                     & {38}                                              & {45}     & {45}     & {45}     \\
                8450                              & 578                                 & 578                                     & {40}                                              & {47}     & {49}     & {48}     \\
                33282                             & 2178                                & 2178                                    & {41}                                              & {49}     & {51}     & {50}     \\
                132098                            & 8450                                & 8450                                    & {41}                                              & {51}     & {52}     & {52}     \\
                526338                            & 33282                               & 33282                                   & {42}                                              & {51}     & {52}     & {54}     \\
                2101250                           & 132098                              & 132098                                  & {43}                                              & {51}     & {54}     & {54}     \\
                8396802                           & 526338                              & 526338                                  & {41}                                              & {51}     & {52}     & {54}     \\
                \bottomrule
            \end{tabular}}
        \caption{Immersed square: MinRes iteration counts for $\bbB_{\mathrm{el}}$. Ratio $h/h_2=0.73$.}
        \label{tab:2d_square_linear_elasticity}
    \end{minipage}\hfill
    \begin{minipage}[t]{0.49\textwidth}
        \vspace{0pt}
        \centering
        \sisetup{group-separator={,}, group-minimum-digits=4}
        \resizebox{\linewidth}{!}{%
            \begin{tabular}{
                    S[table-format=7]
                    S[table-format=7]
                    S[table-format=6]
                    S[table-format=2]
                    S[table-format=2]
                    S[table-format=2]
                    S[table-format=2]
                }
                \toprule
                \multicolumn{7}{c}{\textbf{Elasticity test case (3D, immersed cube)}}                                                                                                                                  \\
                \midrule
                {$\mathrm{dof}_{\boldsymbol{u}}$} & {$\mathrm{dof}_{\boldsymbol{u_2}}$} & {$\mathrm{dof}_{\boldsymbol{\lambda}}$} & \multicolumn{4}{c}{MinRes iterations}                                              \\
                                                  &                                     &                                         & \multicolumn{4}{c}{$\mu_2/\mu=\lambda_2/\lambda$}                                  \\
                \cmidrule(lr){4-7}
                                                  &                                     &                                         & {$10$}                                            & {$10^2$} & {$10^3$} & {$10^4$} \\
                \midrule
                2187                              & 81                                  & 81                                      & {41}                                              & {56}     & {59}     & {61}     \\
                14739                             & 375                                 & 375                                     & {46}                                              & {59}     & {61}     & {61}     \\
                107811                            & 2187                                & 2187                                    & {48}                                              & {64}     & {67}     & {68}     \\
                823875                            & 14739                               & 14739                                   & {46}                                              & {58}     & {62}     & {62}     \\
                6440067                           & 107811                              & 107811                                  & {48}                                              & {62}     & {64}     & {66}     \\
                \bottomrule
            \end{tabular}}
        \caption{Immersed cube: MinRes iteration counts for $\bbBt_{\mathrm{el}}$. Ratio $h/h_2=0.73$.}
        \label{tab:3d_cube_linear_elasticity}
    \end{minipage}
\end{table}

\begin{figure}[ht]
    \centering
    \captionsetup{skip=4pt}
    \captionsetup[subfigure]{skip=3pt,belowskip=0pt}
    \begin{minipage}[t]{0.42\textwidth}
        \centering
        \subcaptionbox{Two-dimensional experiment.\label{fig:displacement_square}}[\linewidth]{%
            \includegraphics[width=0.70\linewidth]{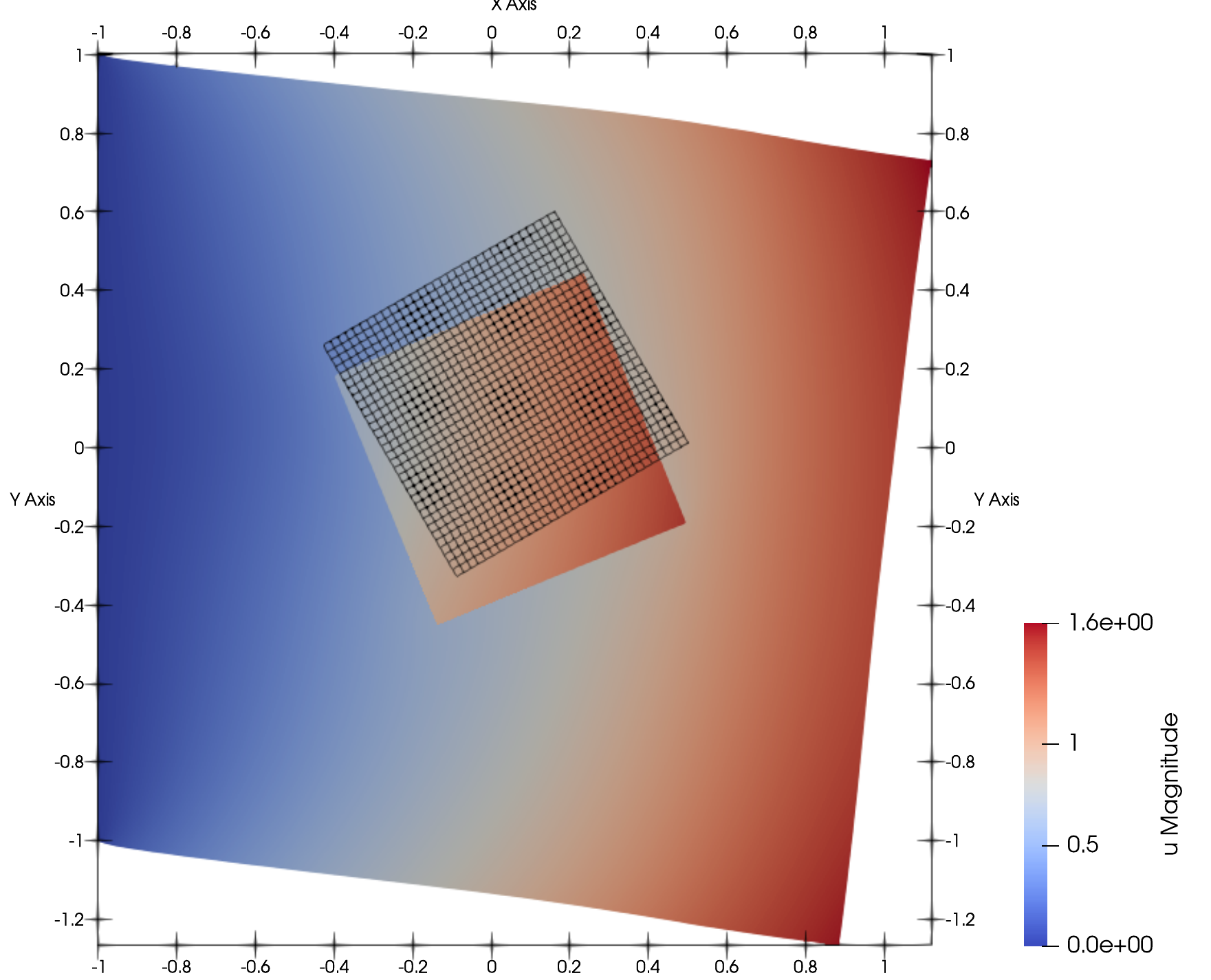}}
    \end{minipage}\hfill
    \begin{minipage}[t]{0.42\textwidth}
        \centering
        \subcaptionbox{Three-dimensional experiment.\label{fig:displacement_cube}}[\linewidth]{%
            \includegraphics[width=0.70\linewidth]{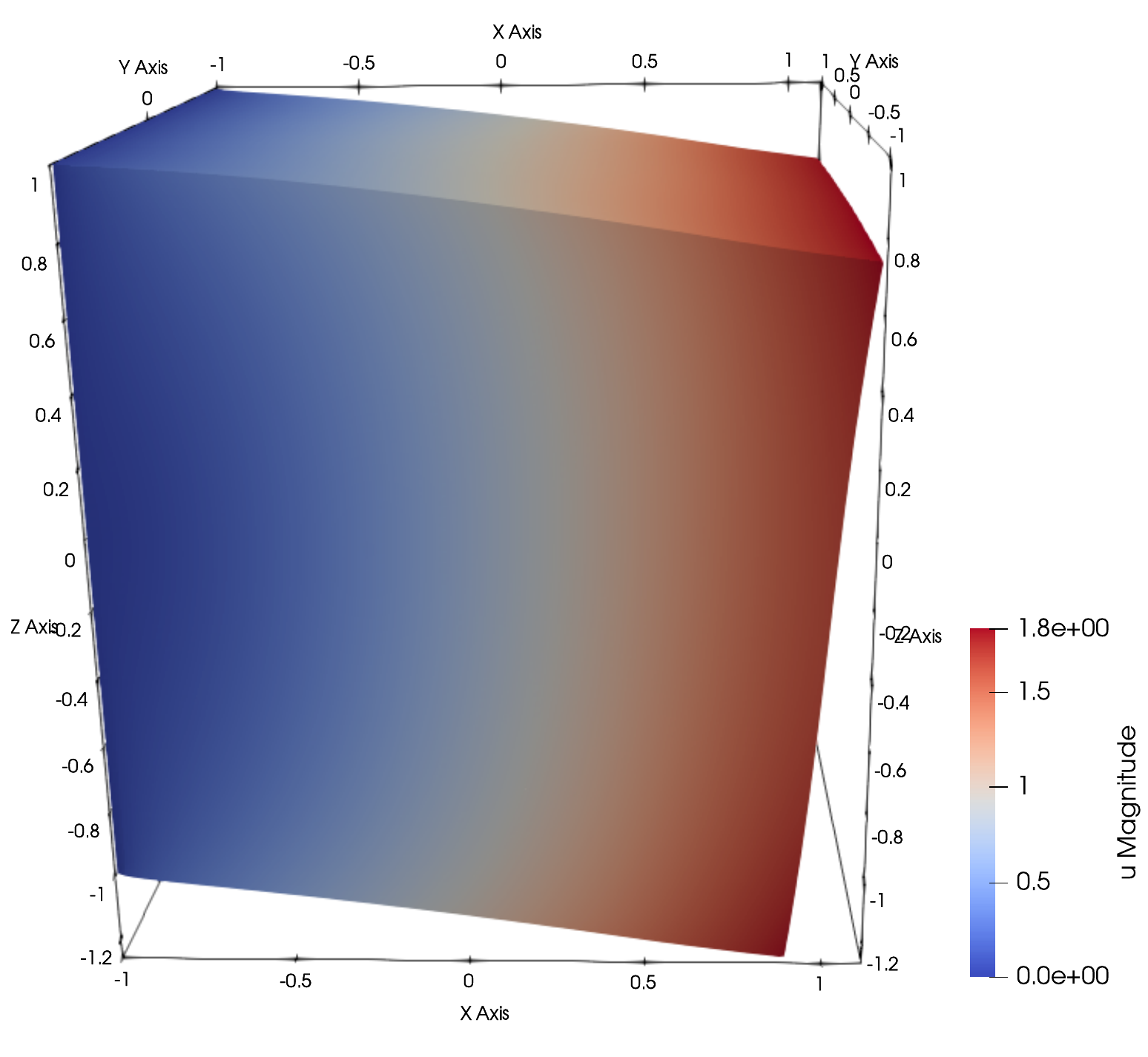}}
    \end{minipage}
    \caption{Magnitude of the displacement field $\boldsymbol{u}_h$ on the background mesh
        when $\mu_2/\mu=\lambda_2/\lambda=10^2$. The immersed mesh in the two-dimensional experiment is shown in wireframe (not shown in three dimensions).}
\end{figure}
\endgroup

\subsection{Large-scale experiment}

We conclude with a strong scaling study on large-scale three-dimensional elasticity problems. To accelerate the time-to-solution, robust iteration counts
can be complemented by efficient operator evaluations and preconditioner applications. For this reason, matrix-free techniques provide a natural way
to address this computational challenge while avoiding the storage of assembled matrices. Indeed, within Krylov methods, sparse matrix--vector products are typically memory-bandwidth limited.
Matrix-free methods avoid storing the matrix by evaluating finite element integrals on the fly, often improving performance for polynomial degrees $p \geq 2$.
For tensor-product bases and quadrature rules, sum-factorization~\cite{ORSZAG198070} reduces the element-wise cost
from $\mathcal{O}((p+1)^{2d})$ to $\mathcal{O}(d(p+1)^{d+1})$ in $d$ spatial dimensions.
Mesh partitioning and SIMD vectorization over batches of elements provide further parallelism~\cite{KRONBICHLER2012135}.

Therefore, we exploit state-of-the-art matrix-free kernels to perform the operator evaluation
of the blocks of the system and of the preconditioner (except for the coupling matrix $\mat{C}$, which is instead explicitly assembled in a setup phase). Notably, the
first two diagonal blocks of $\bbBt_{\mathrm{el}}$ involve the inversion of standard elliptic operators defined on the background and immersed domains, respectively, and can be efficiently inverted using geometric
multigrid methods, which are well-suited for matrix-free implementations.
More precisely, we exploit the geometric multigrid infrastructure of the \textsc{deal.II} library described in~\cite{global_coarsening_mg} to approximately
invert the diagonal blocks of the preconditioner $\bbBt_{\mathrm{el}}$ in a matrix-free fashion. The first two diagonal blocks of $\bbBt_{\mathrm{el}}$ are inverted using 2 V-cycles of geometric multigrid. The smoother is a Chebyshev
accelerated Jacobi scheme~\cite{ADAMS2003593} of degree $6$, which uses the inverse of the matrix diagonal on each
level (precomputed before solving) and the level operator application for computing residuals in the
Jacobi-type iteration. Smoother parameters are determined by an eigenvalue estimation based on Lanczos iterations. The inversion of the mass matrices in the third
block is performed using Chebyshev iteration of fixed degree $10$, with fixed spectral bounds enclosing the whole spectrum. Together with the fixed SPD multigrid approximations for the first two
blocks, this yields a fixed symmetric positive-definite preconditioner suitable for MinRes.

The simulations have been run on the MareNostrum 5 cluster at Barcelona Supercomputing Center, where each compute node features 2x Intel Sapphire Rapids 8480Y+
processors at 2GHz (56 cores each, 112 cores per node) and 256GB of main memory\footnote{Technical characteristics are available at \url{https://www.bsc.es/marenostrum/marenostrum-5}.}. All nodes
have AVX-512 ISA extension, which allows 8 double-precision or 16 single-precision
numbers to be processed per instruction.
We keep the problem size fixed while increasing the number of MPI processes up to $7,168$ on $64$ nodes.

\begin{figure}[!t]
    \centering
    \begin{minipage}[t]{0.49\textwidth}
        \centering
        \subcaptionbox{$\mathcal{Q}^1$--$\mathcal{Q}^1$--$\mathcal{Q}^1$ element.\label{fig:strong-scaling-p1}}[\linewidth]{%
            \includegraphics[width=\linewidth]{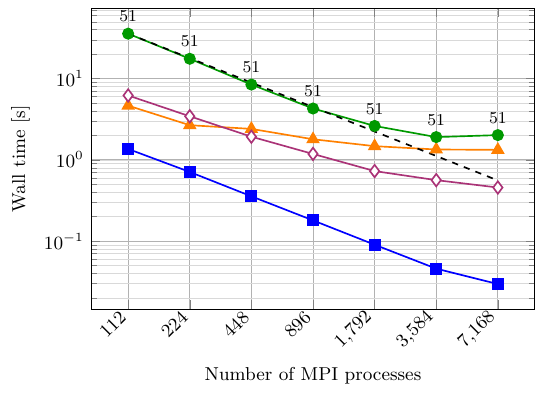}}
    \end{minipage}\hfill
    \begin{minipage}[t]{0.49\textwidth}
        \centering
        \subcaptionbox{$\mathcal{Q}^2$--$\mathcal{Q}^2$--$\mathcal{Q}^1$ element.\label{fig:strong-scaling-p2}}[\linewidth]{%
            \includegraphics[width=\linewidth]{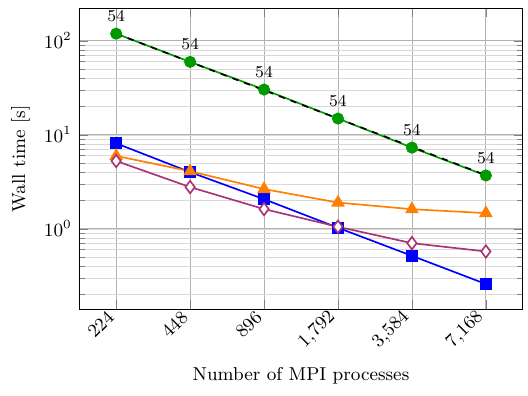}}
    \end{minipage}
    \par\smallskip
    \includegraphics[width=0.85\linewidth]{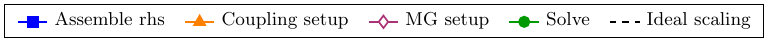}
    \caption{Strong scaling of the solver phases for the three-dimensional linear elasticity problem. Left: $\mathcal{Q}^1$--$\mathcal{Q}^1$--$\mathcal{Q}^1$ element, with
        $51,031,590$ total degrees of freedom. Right: $\mathcal{Q}^2$--$\mathcal{Q}^2$--$\mathcal{Q}^1$ element, with
        $405,948,777$ total degrees of freedom. Iteration counts are reported on top of each data point in the
        solve phase. The dashed line denotes ideal scaling anchored to the first available solve data point. Ratio $h/h_2=0.73$, and $\mu_2/\mu=\lambda_2/\lambda=10^2$.}
    \label{fig:strong-scaling}
\end{figure}

We measure the wall time (in seconds) of the main phases of the solver as a function of the number of MPI processes: assembly of the right-hand side, coupling setup,
multigrid (MG) setup, and solve. Figure~\ref{fig:strong-scaling-p1} reports the strong scaling study for the $\mathcal{Q}^1$--$\mathcal{Q}^1$--$\mathcal{Q}^1$ element, with a fixed problem
size of $51,031,590$ total degrees of freedom. Iteration counts
are reported on top of each data point in the solve phase, and are constant for each number of MPI processes involved. The preconditioner achieves ideal strong scaling, until communication overhead starts to dominate at large core counts. Similarly, ideal scaling is observed for the assembly of the right-hand side. The coupling phase, which
involves the setup and assembly of the coupling term, is less scalable in these tests. Notice that, for moving immersed domains, the coupling matrix typically requires updating at each time step.

To test the scalability of the preconditioner at higher orders, which benefit from matrix-free techniques, we repeat the same experiment for the $\mathcal{Q}^2$--$\mathcal{Q}^2$--$\mathcal{Q}^1$ element\footnote{The previous analysis does not directly cover this element. Its stability in standard norms can be established by adapting the arguments of~\cite{AURICCHIO201536}.}, featuring higher order spaces for $V_h$ and $V_{2,h}$. Since $V_{2,h} \ne \Lambda_h$, the blocks $(2, 3)$ and $(3, 2)$ of $\bbA$ are now rectangular, and the mass matrix of the second block in~\eqref{eq:precond_el} is defined
on $V_{2,h}$, while matrices on the third block of~\eqref{eq:precond_el} are still defined on $\Lambda_h$. The problem's size is fixed at $405,948,777$ degrees of freedom, and results
are reported in Figure~\ref{fig:strong-scaling-p2}. The preconditioner achieves ideal strong scalability from $2$ to $64$ nodes, with constant iteration counts. The same trends as in the equal-order case with $\mathcal{Q}^1$--$\mathcal{Q}^1$--$\mathcal{Q}^1$ element are observed for the assembly of the right-hand side and of the coupling operator.

\section{Conclusions}\label{sec:conclusions}

Using the operator preconditioning framework, we have derived a parameter-robust preconditioner for the fictitious domain discretization with distributed Lagrange
multipliers of the elliptic interface problem, where the multiplier space is defined on an immersed domain. We have analyzed
its robustness with respect to model and discretization parameters, under a mesh ratio condition at the discrete level. Numerical experiments support the theoretical findings and demonstrate lower setup and
solution times than the augmented Lagrangian preconditioner of~\cite{BENZI2026522}, and the scalability
of the preconditioner has been demonstrated through parallel large-scale runs exploiting highly efficient matrix-free techniques.
Ongoing work focuses on extending the present approach to coupled multiphysics problems, in particular to monolithic discretizations of fluid--structure interaction problems using fictitious domain formulations, and on porting
this strategy to GPU architectures.

\section*{Data availability}
The source code used to generate the results presented in this paper is available at the maintained GitHub repository \href{https://github.com/fdrmrc/OP4FDD}{github.com/fdrmrc/OP4FDD}
provided by the author. Detailed instructions for compilation and execution are provided in the repository.

\section*{Acknowledgments}
The author is grateful to Luca Heltai (University of Pisa) for helpful discussions.

\section*{Declarations}
AI tools were used to assist with text editing for spelling and grammar (ChatGPT 6.1 Sol), and to polish the
early version of the code for reproducibility (Claude Opus 5.5). The author takes full responsibility for all content.

\appendix
\section{Auxiliary results}

\begin{lemma}\label{lem:s_bound}
    Let $\theta \coloneqq(1/\beta_1 + 1/\alpha)^{-1}$. Then, for all $\eta\ge 0$,
    \begin{align*}
        \quad & \tfrac12(\theta\eta+\beta_1) \le s(\eta)^{-1} \le \theta\eta+\beta_1,
              &                                                                       & s(\eta)=\frac{1}{\beta_1(\eta+1)}+\frac{1}{\alpha\eta+\beta_1}.
    \end{align*}
\end{lemma}
\begin{proof}
    By definition, $\theta(\beta_1+\alpha)=\beta_1\alpha$, and
    $\theta<\min\{\beta_1,\alpha\}$ since $\alpha/(\beta_1+\alpha)<1$ and
    $\beta_1/(\beta_1+\alpha)<1$.

    For $a,b>0$, $(1/a+1/b)^{-1}=ab/(a+b)$.
    With $a=\beta_1(\eta+1)$ and $b=\alpha\eta+\beta_1$, the identity
    $\theta(\beta_1+\alpha)=\beta_1\alpha$ gives
    $$(\theta\eta+\beta_1)(a+b)-ab=\beta_1(2\theta\eta+\beta_1)\ge 0,$$
    and $$2ab-(\theta\eta+\beta_1)(a+b)=\beta_1\eta\,(\alpha\eta+\beta_1+\alpha-2\theta)\ge 0,$$
    since $\beta_1+\alpha -2\theta = \frac{\betaone^2 + \gap^2}{\betaone + \gap}>0$. Dividing by $a+b>0$ yields the claim.
\end{proof}

\bibliographystyle{siamplain}
\bibliography{refs}

@inproceedings{haznics,
  title     = {Rational approximation preconditioners for multiphysics problems},
  author    = {Budi{\v{s}}a, Ana and Hu, Xiaozhe and Kuchta, Miroslav and Mardal, Kent-Andr{\'e} and Zikatanov, Ludmil},
  booktitle = {Numerical Methods and Applications},
  series    = {Lecture Notes in Computer Science},
  volume    = {13858},
  pages     = {100--113},
  year      = {2023},
  publisher = {Springer Nature Switzerland},
  doi       = {10.1007/978-3-031-32412-3_9}
}

@article{prec1d2d,
  author  = {Kuchta, Miroslav and Nordaas, Magne and Verschaeve, Joris C. G. and Mortensen, Mikael and Mardal, Kent-Andre},
  title   = {Preconditioners for Saddle Point Systems with Trace Constraints Coupling 2{D} and 1{D} Domains},
  journal = {SIAM Journal on Scientific Computing},
  volume  = {38},
  number  = {6},
  pages   = {B962-B987},
  year    = {2016},
  doi     = {10.1137/15M1052822}
}

@book{BrezziFortin,
  author    = {Brezzi, Franco and Fortin, Michel},
  year      = {1991},
  title     = {Mixed and Hybrid Finite Element Methods},
  isbn      = {0387975829},
  series    = {Springer Series in Computational Mathematics},
  volume    = {15},
  publisher = {Springer-Verlag},
  address   = {New York},
  doi       = {10.1007/978-1-4612-3172-1}
}

@article{ALprec,
  author  = { Benzi, Michele and  Olshanskii, Maxim A.},
  title   = {An Augmented {L}agrangian-Based Approach to the {O}seen Problem},
  journal = {SIAM Journal on Scientific Computing},
  volume  = {28},
  number  = {6},
  pages   = {2095-2113},
  year    = {2006},
  doi     = {10.1137/050646421},
  url     = {https://doi.org/10.1137/050646421}
}

@article{BoffiGastaldi2024,
  author  = {Alshehri, Najwa and Boffi, Daniele and Gastaldi, Lucia},
  title   = {Unfitted mixed finite element methods for elliptic interface problems},
  journal = {Numerical Methods for Partial Differential Equations},
  volume  = {40},
  number  = {1},
  pages   = {e23063},
  doi     = {10.1002/num.23063},
  url     = {https://onlinelibrary.wiley.com/doi/abs/10.1002/num.23063},
  year    = {2024}
}

@article{Glowinski,
  title   = {Error analysis of a fictitious domain method applied to a {Dirichlet} problem},
  volume  = {12},
  issn    = {1868-937X},
  url     = {https://doi.org/10.1007/BF03167240},
  doi     = {10.1007/BF03167240},
  number  = {3},
  journal = {Japan Journal of Industrial and Applied Mathematics},
  author  = {Girault, V. and Glowinski, R.},
  year    = {1995},
  pages   = {487--514}
}

@article{BG_numermath,
  author     = {Boffi, Daniele and Gastaldi, Lucia},
  title      = {A fictitious domain approach with {L}agrange multiplier for
                fluid-structure interactions},
  journal    = {Numer. Math.},
  fjournal   = {Numerische Mathematik},
  volume     = {135},
  year       = {2017},
  number     = {3},
  pages      = {711--732},
  issn       = {0029-599X},
  mrclass    = {65M60 (65M12 74F10)},
  mrnumber   = {3606460},
  mrreviewer = {Huazhong Tang},
  doi        = {10.1007/s00211-016-0814-1},
  url        = {https://doi.org/10.1007/s00211-016-0814-1}
}

@article{Trilinos,
  author     = {Heroux, Michael A. and Bartlett, Roscoe A. and Howle, Vicki E. and Hoekstra, Robert J. and Hu, Jonathan J. and Kolda, Tamara G. and Lehoucq, Richard B. and Long, Kevin R. and Pawlowski, Roger P. and Phipps, Eric T. and Salinger, Andrew G. and Thornquist, Heidi K. and Tuminaro, Ray S. and Willenbring, James M. and Williams, Alan and Stanley, Kendall S.},
  title      = {An Overview of the {T}rilinos Project},
  year       = {2005},
  issue_date = {September 2005},
  publisher  = {Association for Computing Machinery},
  journal    = {ACM Transactions on Mathematical Software},
  address    = {New York, NY, USA},
  volume     = {31},
  number     = {3},
  issn       = {0098-3500},
  url        = {https://doi.org/10.1145/1089014.1089021},
  doi        = {10.1145/1089014.1089021},
  pages      = {397-423},
  numpages   = {27}
}

@article{dealIIdesign,
  author  = {Daniel Arndt and Wolfgang Bangerth and Denis Davydov and Timo Heister and Luca Heltai and Martin Kronbichler and Matthias Maier and Jean-Paul Pelteret and Bruno Turcksin and David Wells},
  journal = {Computers {\&} Mathematics with Applications},
  pages   = {407--422},
  title   = {The deal.{II} finite element library: Design, features, and insights},
  volume  = {81},
  year    = {2021},
  doi     = {10.1016/j.camwa.2020.02.022}
}

@article{boffi2023comparison,
  title     = {A comparison of non-matching techniques for the finite element approximation of interface problems},
  author    = {Boffi, Daniele and Cangiani, Andrea and Feder, Marco and Gastaldi, Lucia and Heltai, Luca},
  journal   = {Computers {\&} Mathematics with Applications},
  volume    = {151},
  pages     = {101--115},
  year      = {2023},
  publisher = {Elsevier},
  doi       = {10.1016/j.camwa.2023.09.017}
}

@article{BOFFI2024406,
  title    = {{A parallel solver for fluid-structure interaction problems with {L}agrange multiplier}},
  journal  = {Mathematics and Computers in Simulation},
  volume   = {220},
  pages    = {406-424},
  year     = {2024},
  issn     = {0378-4754},
  doi      = {10.1016/j.matcom.2024.01.027},
  url      = {https://www.sciencedirect.com/science/article/pii/S0378475424000399},
  author   = {Daniele Boffi and Fabio Credali and Lucia Gastaldi and Simone Scacchi}
}

@article{BOFFI2022115650,
  title   = {On the interface matrix for fluid-structure interaction problems with fictitious domain approach},
  journal = {Computer Methods in Applied Mechanics and Engineering},
  volume  = {401},
  pages   = {115650},
  year    = {2022},
  issn    = {0045-7825},
  doi     = {10.1016/j.cma.2022.115650},
  url     = {https://www.sciencedirect.com/science/article/pii/S0045782522006053},
  author  = {Daniele Boffi and Fabio Credali and Lucia Gastaldi}
}

@article{boffi2024quadratureerrorestimatesnonmatching,
  title        = {Quadrature error estimates on non-matching grids in a fictitious domain framework for fluid-structure interaction problems},
  author       = {Daniele Boffi and Fabio Credali and Lucia Gastaldi},
  journal      = {Numerische Mathematik},
  volume       = {158},
  number       = {5},
  pages        = {1895--1932},
  year         = {2026},
  doi          = {10.1007/s00211-026-01553-3},
  primaryclass = {math.NA}
}

@book{boffi2013mixed,
  title     = {Mixed Finite Element Methods and Applications},
  author    = {Boffi, D. and Brezzi, F. and Fortin, M.},
  isbn      = {9783642365195},
  lccn      = {2013940257},
  series    = {Springer Series in Computational Mathematics},
  volume    = {44},
  url       = {https://books.google.it/books?id=mRhAAAAAQBAJ},
  year      = {2013},
  publisher = {Springer Berlin Heidelberg},
  doi       = {10.1007/978-3-642-36519-5}
}

@inproceedings{alshehri2025multigridpreconditioningfddlmmethod,
  title        = {{Multigrid Preconditioning for FD-DLM Method in Elliptic Interface Problems}},
  booktitle    = {Large-Scale Scientific Computations},
  series       = {Lecture Notes in Computer Science},
  pages        = {103--111},
  publisher    = {Springer Nature Switzerland},
  author       = {Najwa Alshehri and Daniele Boffi and Chayapol Chaoveeraprasit},
  year         = {2026},
  doi          = {10.1007/978-3-032-22221-3_11},
  primaryclass = {math.NA}
}

@article{GMRES_eigs,
  author  = {Greenbaum, Anne and Pt\'{a}k, Vlastimil and Strako\v{s}, Zden\v{e}k},
  title   = {Any Nonincreasing Convergence Curve is Possible for {GMRES}},
  journal = {SIAM Journal on Matrix Analysis and Applications},
  volume  = {17},
  number  = {3},
  pages   = {465-469},
  year    = {1996},
  doi     = {10.1137/S0895479894275030},
  url     = {https://doi.org/10.1137/S0895479894275030}
}

@article{Peskin_2002,
  title   = {The immersed boundary method},
  volume  = {11},
  doi     = {10.1017/S0962492902000077},
  journal = {Acta Numerica},
  author  = {Peskin, Charles S.},
  year    = {2002},
  pages   = {479-517}
}

@article{Farrell2019,
  author  = {Farrell, Patrick E. and Mitchell, Lawrence and Wechsung, Florian},
  title   = {An Augmented {L}agrangian Preconditioner for the {3D} Stationary Incompressible {N}avier-{S}tokes Equations at High {R}eynolds Number},
  journal = {SIAM Journal on Scientific Computing},
  volume  = {41},
  number  = {5},
  pages   = {A3073-A3096},
  year    = {2019},
  doi     = {10.1137/18M1219370},
  url     = {https://doi.org/10.1137/18M1219370}
}

@book{Ern_Guermond-FEM-2004,
  author    = {Ern, A. and Guermond, J.-L.},
  title     = {{Theory and Practice of Finite Elements}},
  series    = {Applied Mathematical Sciences},
  volume    = {159},
  publisher = {Springer-Verlag},
  address   = {New York},
  year      = {2004},
  isbn      = {9781475743555},
  doi       = {10.1007/978-1-4757-4355-5}
}

@article{AURICCHIO201536,
  title    = {On a fictitious domain method with distributed {L}agrange multiplier for interface problems},
  journal  = {Applied Numerical Mathematics},
  volume   = {95},
  pages    = {36-50},
  year     = {2015},
  issn     = {0168-9274},
  doi      = {10.1016/j.apnum.2014.05.013},
  url      = {https://www.sciencedirect.com/science/article/pii/S0168927414001238},
  author   = {Ferdinando Auricchio and Daniele Boffi and Lucia Gastaldi and Adrien Lefieux and Alessandro Reali}
}

@article{BOFFI20142151,
  title   = {Mixed formulation for interface problems with distributed {L}agrange multiplier},
  journal = {Computers {\&} Mathematics with Applications},
  volume  = {68},
  number  = {12, Part B},
  pages   = {2151-2166},
  year    = {2014},
  note    = {Advances in Computational Partial Differential Equations},
  issn    = {0898-1221},
  doi     = {10.1016/j.camwa.2014.07.020},
  author  = {Daniele Boffi and Lucia Gastaldi and Michele Ruggeri}
}

@article{REGAZZONI2024117327,
  title    = {An optimally convergent Fictitious Domain method for interface problems},
  journal  = {Computer Methods in Applied Mechanics and Engineering},
  volume   = {431},
  pages    = {117327},
  year     = {2024},
  issn     = {0045-7825},
  doi      = {10.1016/j.cma.2024.117327},
  url      = {https://www.sciencedirect.com/science/article/pii/S0045782524005826},
  author   = {Francesco Regazzoni}
}

@article{BENZI2026118522,
  title    = {Scalable augmented {L}agrangian preconditioners for fictitious domain problems},
  journal  = {Computer Methods in Applied Mechanics and Engineering},
  volume   = {450},
  pages    = {118522},
  year     = {2026},
  issn     = {0045-7825},
  doi      = {10.1016/j.cma.2025.118522},
  url      = {https://www.sciencedirect.com/science/article/pii/S0045782525007947},
  author   = {Michele Benzi and Marco Feder and Luca Heltai and Federica Mugnaioni}
}

@article{BoffiCavalliniGastaldi2015,
  author  = {Boffi, Daniele and Cavallini, Nicola and Gastaldi, Lucia},
  title   = {The Finite Element Immersed Boundary Method with Distributed {L}agrange Multiplier},
  journal = {SIAM Journal on Numerical Analysis},
  volume  = {53},
  number  = {6},
  pages   = {2584-2604},
  year    = {2015},
  doi     = {10.1137/140978399},
  url     = {https://doi.org/10.1137/140978399}
}

@article{BoffiConditioning,
  author  = {Boffi, Daniele and Credali, Fabio and Gastaldi, Lucia},
  title   = {On the Stability and Conditioning of a Fictitious Domain Formulation for Fluid-Structure Interaction Problems},
  journal = {SIAM Journal on Scientific Computing},
  volume  = {48},
  number  = {4},
  pages   = {A1843-A1869},
  year    = {2026},
  doi     = {10.1137/25M1758817},
  url     = {https://doi.org/10.1137/25M1758817}
}

@article{schoeberl_1999,
  author  = {Sch\"{o}berl, Joachim},
  year    = {1999},
  month   = {11},
  pages   = {97-119},
  title   = {Multigrid Methods for a Parameter Dependent Problem in Primal Variables},
  volume  = {84},
  journal = {Numerische Mathematik},
  doi     = {10.1007/s002110050465}
}

@article{ALZETTA2020106334,
  title   = {Multiscale modeling of fiber reinforced materials via non-matching immersed methods},
  journal = {Computers \& Structures},
  volume  = {239},
  pages   = {106334},
  year    = {2020},
  issn    = {0045-7949},
  doi     = {10.1016/j.compstruc.2020.106334},
  author  = {Giovanni Alzetta and Luca Heltai}
}

@article{BENZI2026522,
  title   = {Augmented {L}agrangian preconditioners for fictitious domain formulations of elliptic interface problems},
  journal = {Computers \& Mathematics with Applications},
  volume  = {224},
  pages   = {522-551},
  year    = {2026},
  issn    = {0898-1221},
  doi     = {10.1016/j.camwa.2026.10.042},
  author  = {Michele Benzi and Marco Feder and Luca Heltai and Federica Mugnaioni}
}

@article{FIRMBACH2027119425,
  title    = {Modified augmented {Lagrangian} preconditioning for mixed-dimensional beam-solid coupling},
  journal  = {Computer Methods in Applied Mechanics and Engineering},
  volume   = {463},
  pages    = {119425},
  year     = {2027},
  issn     = {0045-7825},
  doi      = {10.1016/j.cma.2026.119425},
  author   = {Max Firmbach and Ivo Steinbrecher and Alexander Popp and Matthias Mayr}
}

@article{MW,
  author  = {Mardal, Kent-Andre and Winther, Ragnar},
  title   = {Preconditioning discretizations of systems of partial differential equations},
  journal = {Numerical Linear Algebra with Applications},
  volume  = {18},
  number  = {1},
  pages   = {1--40},
  year    = {2011},
  doi     = {10.1002/nla.716}
}

@article{MUMPS,
  title   = {{Performance and Scalability of the Block Low-Rank Multifrontal Factorization on Multicore Architectures}},
  author  = {P.R. Amestoy and A. Buttari and J.-Y. L'Excellent and T. Mary},
  journal = {ACM Transactions on Mathematical Software},
  volume  = 45,
  number  = {1},
  pages   = {2:1--2:26},
  year    = {2019},
  doi     = {10.1145/3242094}
}

@article{ADAMS2003593,
  title    = {Parallel multigrid smoothing: polynomial versus {G}auss-{S}eidel},
  journal  = {Journal of Computational Physics},
  volume   = {188},
  number   = {2},
  pages    = {593-610},
  year     = {2003},
  issn     = {0021-9991},
  doi      = {10.1016/S0021-9991(03)00194-3},
  url      = {https://www.sciencedirect.com/science/article/pii/S0021999103001943},
  author   = {Mark Adams and Marian Brezina and Jonathan Hu and Ray Tuminaro}
}

@article{KRONBICHLER2012135,
  title   = {A generic interface for parallel cell-based finite element operator application},
  journal = {Computers {\&} Fluids},
  volume  = {63},
  pages   = {135-147},
  year    = {2012},
  issn    = {0045-7930},
  doi     = {10.1016/j.compfluid.2012.04.012},
  url     = {https://www.sciencedirect.com/science/article/pii/S0045793012001429},
  author  = {Martin Kronbichler and Katharina Kormann}
}

@article{ORSZAG198070,
  title   = {Spectral methods for problems in complex geometries},
  journal = {Journal of Computational Physics},
  volume  = {37},
  number  = {1},
  pages   = {70-92},
  year    = {1980},
  issn    = {0021-9991},
  doi     = {10.1016/0021-9991(80)90005-4},
  url     = {https://www.sciencedirect.com/science/article/pii/0021999180900054},
  author  = {Steven A Orszag}
}

@article{global_coarsening_mg,
  author     = {Munch, Peter and Heister, Timo and Prieto Saavedra, Laura and Kronbichler, Martin},
  title      = {Efficient Distributed Matrix-free Multigrid Methods on Locally Refined Meshes for {FEM} Computations},
  year       = {2023},
  issue_date = {March 2023},
  publisher  = {Association for Computing Machinery},
  address    = {New York, NY, USA},
  volume     = {10},
  number     = {1},
  issn       = {2329-4949},
  url        = {https://doi.org/10.1145/3580314},
  doi        = {10.1145/3580314},
  journal    = {ACM Trans. Parallel Comput.},
  month      = mar,
  articleno  = {3},
  numpages   = {38}
}

@article{PaigeSaunders1975,
  author  = {Paige, C. C. and Saunders, M. A.},
  title   = {Solution of Sparse Indefinite Systems of Linear Equations},
  journal = {SIAM Journal on Numerical Analysis},
  volume  = {12},
  number  = {4},
  pages   = {617-629},
  year    = {1975},
  doi     = {10.1137/0712047},
  url     = {https://doi.org/10.1137/0712047}
}

@book{grisvard2011elliptic,
  title     = {Elliptic Problems in Nonsmooth Domains},
  author    = {Grisvard, P.},
  isbn      = {9781611972023},
  lccn      = {2011028564},
  series    = {Classics in Applied Mathematics},
  volume    = {69},
  url       = {https://books.google.it/books?id=lJTtujsZSgkC},
  year      = {2011},
  publisher = {Society for Industrial and Applied Mathematics},
  address   = {Philadelphia},
  doi       = {10.1137/1.9781611972030}
}

@book{nicaise,
  author    = {Nicaise, Serge},
  title     = {Polygonal Interface Problems},
  volume    = {39},
  series    = {Methoden und Verfahren der mathematischen Physik},
  publisher = {Peter Lang},
  address   = {Frankfurt am Main},
  year      = {1993},
  isbn      = {9783631463802},
  url       = {https://catalogue.bnf.fr/ark:/12148/cb391288606}
}

@article{CrouzeixThomee1987,
  issn      = {00255718, 10886842},
  url       = {http://www.jstor.org/stable/2007825},
  author    = {M. Crouzeix and V. Thomée},
  journal   = {Mathematics of Computation},
  number    = {178},
  pages     = {521--532},
  publisher = {American Mathematical Society},
  title     = {The Stability in {$L_p$} and {$W_p^1$} of the {$L_2$}-Projection onto Finite Element Function Spaces},
  volume    = {48},
  year      = {1987},
  doi       = {10.1090/S0025-5718-1987-0878688-2}
}

@article{hu2026parameter,
  title        = {Parameter-robust preconditioner for {Stokes-Darcy} coupled problem with {Lagrange} multiplier},
  author       = {Xiaozhe Hu and Miroslav Kuchta and Kent-Andre Mardal and Xue Wang},
  year         = {2026},
  journal      = {SIAM Journal on Scientific Computing},
  volume       = {48},
  number       = {2},
  pages        = {A878--A904},
  doi          = {10.1137/25M1744393},
  primaryclass = {math.NA}
}

@article{cai2026parameter,
  title        = {Parameter-Robust Preconditioners for A Four-Field Thermo-Poroelasticity Model},
  author       = {Mingchao Cai and Miroslav Kuchta and Jingzhi Li and Ziliang Li and Kent-Andre Mardal},
  year         = {2026},
  journal      = {SIAM Journal on Scientific Computing},
  volume       = {48},
  number       = {1},
  pages        = {A49--A73},
  doi          = {10.1137/25M1726923},
  primaryclass = {math.NA}
}

@article{sande2025robust,
  title         = {On a robust inf-sup condition for the {Stokes} problem in slender domains--with application to preconditioning},
  author        = {Sande, Espen and Koch, Timo and Kuchta, Miroslav and Mardal, Kent-Andre},
  year          = {2025},
  eprint        = {2510.24590},
  journal       = {arXiv:2510.24590},
  archiveprefix = {arXiv},
  primaryclass  = {math.NA}
}

@incollection{dimola2024robust,
  title     = {Robust Preconditioning of Mixed-Dimensional {PDEs} on {3D-1D} Domains Coupled with {L}agrange Multipliers},
  author    = {Dimola, Nunzio and Kuchta, Miroslav and Mardal, Kent-Andre and Zunino, Paolo},
  booktitle = {Quantitative Approaches to Microcirculation: Mathematical Models, Computational Methods and Data Analysis},
  series    = {SEMA SIMAI Springer Series},
  pages     = {137--171},
  year      = {2024},
  publisher = {Springer Nature Switzerland},
  doi       = {10.1007/978-3-031-58519-7_6}
}

@article{benzi2026unifiedvariationalframeworkaugmented,
  title         = {A unified variational framework for augmented {L}agrangian preconditioning of saddle point systems: applications to mixed-dimensional and fictitious domain methods},
  author        = {Michele Benzi and Marco Feder and Luca Heltai and Federica Mugnaioni},
  journal       = {arXiv:2610.30357},
  year          = {2026},
  eprint        = {2610.30357},
  archiveprefix = {arXiv},
  primaryclass  = {math.NA}
}

@book{evans2022partial,
  title     = {Partial differential equations},
  author    = {Evans, Lawrence C},
  volume    = {19},
  year      = {2022},
  publisher = {American Mathematical Society}
}

\end{document}